\documentclass[11pt]{amsart}
\usepackage[margin=22mm]{geometry}

\usepackage{mathtools}
\usepackage{amsmath}
\usepackage{amssymb}
\usepackage{amsthm}
\usepackage{amsfonts}
\usepackage{enumerate}
\usepackage{url}
\usepackage{xcolor}

\DeclareMathOperator{\Ob}{\mathrm{Ob}}
\DeclareMathOperator{\Ima}{\mathrm{Im}}
\DeclareMathOperator{\Ker}{\mathrm{Ker}}

\newcommand{\Z}{\mathbb{Z}}
\newcommand{\Crs}{\mathrm{Crs}}
\newcommand{\meri}{\mathfrak{m}}
\newcommand{\longi}{\mathfrak{l}}

\theoremstyle{plain}
\newtheorem{Theorem}{Theorem}[section]
\newtheorem{Proposition}[Theorem]{Proposition}
\newtheorem{Lemma}[Theorem]{Lemma}

\theoremstyle{definition}
\newtheorem{Definition}[Theorem]{Definition}
\newtheorem{Remark}[Theorem]{Remark}

\title{Identities among Relations for Rational Dehn Fillings of Link Exteriors}
\author{Koki Yanagida}
\address{Faculty of Science and Engineering, Saga University, 1 Honjo-machi, Saga 840--8502, Japan}
\email{yngdk127@cc.saga-u.ac.jp}
\date{}

\begin{document}
\begin{abstract}
Some invariants of closed $3$-manifolds are computed from the equivariant cellular chain complex of the universal cover.
To determine the boundary map in degree $3$ from a group presentation, one needs an identity among relations in addition to the presentation.
Given data for a distinguished boundary torus of a link exterior and any filling slope, we give a uniform formula for the boundary of the new degree-$3$ generator in the fundamental crossed complex of a cellular model for the Dehn filling.
After passing to a reduced model, this boundary gives an identity among relations for the group presentation induced by the model after filling.
The construction combines a continued-fraction model of the filling solid torus with the higher homotopy Seifert--van Kampen theorem for crossed complexes.
Moreover, by choosing a model of the link exterior and boundary torus data for a surgery presentation, we can construct a group presentation and an identity among relations for every closed connected oriented $3$-manifold.

\end{abstract}

\maketitle


\medskip

\section{Introduction}
\label{sec:introduction}

This work begins with the problem of explicitly computing invariants determined by chain complexes of closed $3$-manifolds.
One method constructs the equivariant cellular chain complex of the universal cover directly from a presentation of the fundamental group.
After fixing the generators and relators, this method determines the boundary maps in degrees $1$ and $2$ by Fox differentiation.
The presentation alone, however, does not determine the boundary map in degree $3$.
The attaching map of each $3$-cell yields an identity among relations that determines the corresponding column of the third boundary map of the universal cover \cite{Tro62,Sie86}.
The resulting cellular chain complex through degree $3$ also provides input for finding a chain representative of the fundamental $3$-class \cite{Nos22}.
Such cellular chain complexes have been used in explicit computations of twisted Blanchfield pairings \cite{MP18,BCP24,Yan26}, twisted and non-acyclic Reidemeister torsions \cite{Wak21,Wak24,NWY25}, and Dijkgraaf--Witten invariants \cite{Nos22}.
Thus, a uniform construction of identities among relations provides common input for computing these invariants from group presentations.

Previous work on identities among relations includes the introduction of crossed modules, the development of the algebraic theory of identities among relations, and applications to $3$-manifolds.
Whitehead introduced crossed modules as an algebraic structure on relative second homotopy groups and described the relative second homotopy group obtained by attaching $2$-cells as a free crossed module \cite{Whi49}.
The algebraic theory of identities among relations began with Peiffer and Reidemeister \cite{Pei49,Rei49}.
Brown--Huebschmann, Collins--Huebschmann, Sieradski, and others developed the theory relating identities among relations to presentation complexes and the topology of $3$-manifolds \cite{BH82,CH82,Sie86}.
For knotted surfaces in $S^4$, Faria Martins computed the attaching classes of the $3$-handles from a knot with bands and used them to present the fundamental crossed module of the complement \cite[Theorem~2.6]{FM09}.
Nosaka explicitly constructed identities among relations for several closed $3$-manifolds \cite{Nos22}.
After treating integral surgery, Nosaka left the explicit construction of a taut identity for $p/q$-surgery, for arbitrary $p/q\in\mathbb Q$, as a problem for future study \cite[p.~99]{Nos22}.

In this paper, we construct a formula for identities among relations that applies to Dehn fillings of link exteriors.
Let $L\subset S^3$ be a link, and let $E_L:=S^3\setminus\operatorname{int}(\nu L)$ be its exterior.
Our main result is a uniform substitution formula for Dehn filling.
Given a boundary torus class of a distinguished boundary component of $E_L$ and a slope $p/q\in\mathbb Q\cup\{1/0\}$, the formula gives the boundary in degree $3$ of the fundamental crossed complex of a filtered CW model for the $p/q$-Dehn filling (Theorem~\ref{thm:uniform-substitution-formula-for-p-over-q-surgery}).
In a reduced model, this boundary gives an identity among relations for the group presentation after filling (Remark~\ref{rem:identities-among-relations}).
Every knot exterior admits a handle decomposition with one $0$-handle and no $3$-handles.
Choose the filtered CW model and the boundary torus class compatibly with such a decomposition.
After filling, the relative model completes the decomposition to a Heegaard decomposition with a unique $3$-handle.
For every $p/q\in\mathbb Q$, our formula computes the class of its attaching sphere in the free crossed module of the resulting presentation (Remark~\ref{rem:identities-among-relations}).
When several boundary components are filled, the relative models of the filling solid tori may be glued simultaneously to the original link exterior by a pushout, and the same formula gives the boundary of the new degree-$3$ generator for each component.
In particular, the Lickorish--Wallace theorem states that every closed connected oriented $3$-manifold is obtained by integral surgery on a link in $S^3$ \cite{Wal60,Lic62}.
Therefore, after choosing a surgery presentation, a filtered CW model of its link exterior, and a boundary torus class for each component, our formula constructs a group presentation of any closed connected oriented $3$-manifold together with an identity among relations for that presentation.

The proof of the main theorem recasts the geometric gluing that represents a Dehn filling as a gluing problem for crossed complexes and then computes the boundary in degree $3$ after filling.
In a $3$-dimensional free crossed complex, the boundary of a degree-$3$ generator can be read as an identity among relations for the group presentation of a reduced model (Remark~\ref{rem:identities-among-relations}).
The higher homotopy Seifert--van Kampen theorem for crossed complexes, on the other hand, computes the fundamental crossed complex of a gluing of suitably filtered CW complexes as a pushout of crossed complexes \cite[Theorem~8.2.5]{BHS11}.
We therefore regard a Dehn filling as the gluing of a link exterior and a filling solid torus along their boundary tori.
On the link-exterior side, we determine a boundary torus class (Definition~\ref{def:boundary-torus-class}).
On the filling-solid-torus side, we construct a CW model relative to the boundary torus (Proposition~\ref{prop:crossed-complex-of-the-filling-solid-torus}).
We then use these data to glue the fundamental crossed complexes of the link exterior and the filling solid torus over the boundary torus.
The resulting pushout gives the boundary in degree $3$ after filling and yields the formula in the main theorem (Theorem~\ref{thm:uniform-substitution-formula-for-p-over-q-surgery}).

The paper is organized as follows.
Section~\ref{sec:crossed-complexes-in-dimension-3} prepares the correspondence between free crossed complexes and handle decompositions.
Section~\ref{sec:boundary-torus-classes-for-a-distinguished-boundary-component} defines boundary torus classes and derives the $p/q$-filling formula from continued fractions.
Section~\ref{sec:torus-gluings-and-operations} computes the crossed complex of a satellite exterior by torus gluing.
Section~\ref{sec:concrete-models-and-computations} proves an assembly formula based on a rectangle decomposition of a boundary annulus and applies it to braid exteriors and, using certain nonorientable spanning surfaces, to knot exteriors.

\section{Handle Decompositions and Free Crossed Complexes}
\label{sec:crossed-complexes-in-dimension-3}

This section clarifies the correspondence between a finite handle decomposition of a compact connected $3$-manifold and a free basis of its fundamental crossed complex.
We first define a $3$-dimensional free crossed complex algebraically and then construct a filtered CW model in which each $i$-handle corresponds to one $i$-cell.
Proposition~\ref{prop:handle-classes-free-basis} shows that each $i$-handle in this model gives a free basis element in degree $i$ of the fundamental crossed complex.
By a filtered CW model of a space $X$ we mean a CW complex $K_X$ with its skeletal filtration and a chosen homotopy equivalence $h_X\colon K_X\to X$.
In a model relative to a common subspace $A$, the map $h_X$ and its chosen homotopy inverse restrict to the identity on $A$, and both homotopies are stationary on $A$.
See \cite{BHS11,BH90} for details.

We use the following notation and conventions throughout the paper.
For a small category $C$, write $\Ob(C)$ for its set of objects and $C(x,y)$ for the set of morphisms from $x$ to $y$.
We write composition in algebraic order: for $c_1\in C(x,y)$ and $c_2\in C(y,z)$, the product $c_1c_2$ denotes the morphism that first follows $c_1$ and then $c_2$.
We assume that every filtered CW complex used in this paper has a connected filtration in the sense of \cite[Definition~8.1.2]{BHS11}.

\subsection{$3$-Dimensional Free Crossed Complexes}
\label{subsec:free-crossed-complexes-in-dimension-3}

A groupoid is a small category in which every morphism is invertible.
For a groupoid $G$, write $G(x,y)$ for the set of morphisms from $x$ to $y$, $G(x):=G(x,x)$ for the vertex group at $x$, and $1_x$ for the identity morphism at $x$.
The groupoid $G$ is called reduced if it has one object and totally disconnected if $G(x,y)=\emptyset$ whenever $x\neq y$.

We define a left action of a groupoid on a totally disconnected groupoid as follows.
Let $G$ and $H$ be groupoids with the same object set, and suppose that $H$ is totally disconnected.
A left action of $G$ on $H$ is a family of group isomorphisms
\[
g\cdot(-)\colon H(y)\longrightarrow H(x)
\qquad
(g\in G(x,y))
\]
satisfying
\[
\begin{aligned}
1_x\cdot h&=h
&&\bigl(h\in H(x)\bigr),\\
(gg')\cdot h&=g\cdot(g'\cdot h)
&&\bigl(g\in G(x,y),\ g'\in G(y,z),\ h\in H(z)\bigr),\\
g\cdot(h_1h_2)&=(g\cdot h_1)(g\cdot h_2)
&&\bigl(g\in G(x,y),\ h_1,h_2\in H(y)\bigr).
\end{aligned}
\]

In this paper, a $3$-dimensional crossed complex means a crossed complex in the sense of \cite[Chapter~7]{BHS11} that is trivial in degrees at least $4$.
It consists of a sequence
\[
C_3\xrightarrow{\partial_3}C_2\xrightarrow{\partial_2}C_1
\]
of groupoids with a common object set $C_0$, together with left actions of $C_1$ on $C_2$ and $C_3$.
Here $C_1$ is a groupoid, whereas $C_2$ and $C_3$ are totally disconnected groupoids, and the following three conditions hold.
The map $\partial_2\colon C_2\to C_1$ and the action of $C_1$ form a crossed module; thus, they satisfy the compatibility condition
\[
\partial_2(g\cdot c)=g\partial_2(c)g^{-1}
\qquad
\bigl(g\in C_1(x,y),\ c\in C_2(y)\bigr)
\]
and the Peiffer identity
\[
\partial_2(c)\cdot c'=cc'c^{-1}
\qquad
(c,c'\in C_2(x)).
\]
Each $C_3(x)$ is abelian, and $\partial_3$ is $C_1$-equivariant and satisfies
\[
\partial_2\partial_3(u)=1_x
\qquad
(u\in C_3(x)).
\]
Finally, $\Ima(\partial_2)$ acts trivially on $C_3$.

The Peiffer identity implies that the kernel of the boundary map in degree $2$ is central.
Indeed, setting $c=k\in\Ker(\partial_2)\cap C_2(x)$ in that identity gives $kc'k^{-1}=c'$, and hence $\Ker(\partial_2)\cap C_2(x)\subseteq Z(C_2(x))$.
Moreover, for $c\in C_2(x)$ and $k\in\Ker(\partial_2)\cap C_2(x)$, we have $\partial_2(c)\cdot k=ckc^{-1}=k$.
Thus, $\Ima(\partial_2)$ acts trivially on $\Ker(\partial_2)$.
This centrality allows us to define a boundary map from a free abelian group in degree $3$ to $C_2$.

We now describe the free $3$-dimensional crossed complex determined by sets $E_1,E_2,E_3$ of generators and by boundary elements for the generators in degrees $2$ and $3$.
Let $E_0$ be an object set and $E_1$ a directed graph on $E_0$, and assign sets $E_2(x)$ and $E_3(x)$ to each $x\in E_0$.
Let $C_1$ be the free groupoid on $(E_0,E_1)$.
For each $u\in E_2(x)$, choose an attaching word $d_2(u)\in C_1(x)$; this word will be the boundary of the generator $u$.
Let $C_2\to C_1$ be the free crossed $C_1$-module generated by $E_2$, with $\partial_2(u)=d_2(u)$ for every $u\in E_2$.
See \cite[Section~2]{Whi49} for the original definition and \cite[Section~7.3.ii, especially Proposition~7.3.7]{BHS11} for the groupoid formulation.

The degree-$3$ part is determined freely by prescribed elements in the kernel of $\partial_2$.
For each $v\in E_3(x)$, specify
\[
d_3(v)\in\Ker(\partial_2)\cap C_2(x),
\]
and let $\pi_1:=C_1/\Ima(\partial_2)$ be the quotient groupoid with object set $E_0$.
Let $C_3$ be the free left $\pi_1$-module with basis $E_3$.
For $v\in E_3(y)$ and $g\in C_1(x,y)$, write $\bar g$ for the image of $g$ in $\pi_1(x,y)$ and define
\[
\partial_3(\bar g\cdot v):=g\cdot d_3(v).
\]
Changing the representative $g$ changes the right-hand side by the action of $\Ima(\partial_2)$, which is trivial on $\Ker(\partial_2)$.
Thus, the right-hand side is independent of the representative.
Each element $g\cdot d_3(v)$ lies in $\Ker(\partial_2)\cap C_2(x)$ and hence in the center of $C_2(x)$.
The preceding assignment therefore determines a unique homomorphism $\partial_3$ from the free abelian vertex group $C_3(x)$ to $C_2(x)$.
We call the resulting crossed complex the free crossed complex on $E_*=(E_0,E_1,E_2,E_3)$ and call $E_*$ its basis.

\subsection{Free Bases from Handle Decompositions}
\label{subsec:fundamental-crossed-complexes-from-handle-decompositions}

We next associate a filtered CW model and a fundamental crossed complex to a handle decomposition.
Let $X$ be a compact connected $3$-manifold, possibly with boundary, and fix a finite handle decomposition of $X$.
Let $B_0$ be the set obtained by choosing one point in each $0$-handle, and let $X_{\leq i}$ be the union of the handles of index at most $i$.
The standard correspondence between handle decompositions and CW complexes allows us to choose a finite CW complex $K_X$ with $0$-skeleton $B_0$ and one $i$-cell for each $i$-handle.
Write
\[
h_X\colon K_X\longrightarrow X
\]
for a filtration-preserving homotopy equivalence from the skeletal filtration of $K_X$ to the handle filtration of $X$.
We use the skeletal filtration of $K_X$ below and write
\[
\Pi X:=\Pi(K_X^{(0)}\subset K_X^{(1)}\subset K_X^{(2)}\subset K_X).
\]
More explicitly, the groupoid in degree $1$ of $\Pi X$ and the groups in degrees $2$ and $3$ at $x\in B_0$ are
\[
(\Pi X)_1=\pi_1(K_X^{(1)},B_0),
\qquad
(\Pi X)_2(x)=\pi_2(K_X^{(2)},K_X^{(1)},x),
\qquad
(\Pi X)_3(x)=\pi_3(K_X,K_X^{(2)},x).
\]
Here $\partial_2$ sends the class of a relative disk to the class of its boundary loop:
\[
\partial_2\colon
\pi_2(K_X^{(2)},K_X^{(1)},x)
\longrightarrow
\pi_1(K_X^{(1)},x).
\]
The map $\partial_3$ is the composite of the relative homotopy boundary map and the map induced by inclusion:
\[
\pi_3(K_X,K_X^{(2)},x)
\longrightarrow
\pi_2(K_X^{(2)},x)
\longrightarrow
\pi_2(K_X^{(2)},K_X^{(1)},x).
\]
The action of $(\Pi X)_1$ is obtained by changing basepoints along paths in $K_X^{(1)}$, a construction called whiskering in the language of crossed complexes; compare \cite[pp.~341--342, 345--346]{Hat02}.
Together with this action, these groups and boundary maps form the fundamental crossed complex $\Pi X$ of this filtration \cite[Section~7.1.v, especially Definition~7.1.27]{BHS11}.

Each handle determines a basis element as follows.
Write $g_a$ for the degree-$1$ homotopy class represented by an oriented core of a $1$-handle.
For each $2$-handle and $3$-handle, choose a characteristic map and a basepoint path (whisker), and write $r_b$ and $s_c$, respectively, for the corresponding homotopy classes.
Then $\partial_2(r_b)$ is the attaching loop of the $2$-handle, and $\partial_3(s_c)$ is the element in degree $2$ represented by the attaching sphere of the $3$-handle.

\begin{Proposition}
\label{prop:handle-classes-free-basis}
The elements $g_a$, $r_b$, and $s_c$ defined from the $1$-, $2$-, and $3$-handles, respectively, freely generate $\Pi X$ in the following sense.
The groupoid $(\Pi X)_1$ is freely generated by the elements $g_a$.
The map $\partial_2\colon(\Pi X)_2\to(\Pi X)_1$ is the free crossed module generated by the elements $r_b$, with $\partial_2(r_b)$ equal to the attaching loop of the corresponding $2$-handle.
The groupoid $(\Pi X)_3$ is the free left $(\Pi X)_1/\Ima(\partial_2)$-module generated by the elements $s_c$, with $\partial_3(s_c)$ equal to the degree-$2$ element represented by the attaching sphere of the corresponding $3$-handle.

\end{Proposition}

\begin{proof}
We apply the freeness theorem for fundamental crossed complexes to the skeletal filtration of $K_X$.
Under the standard correspondence between handle decompositions and CW complexes, each relative handle attachment is replaced by a cell attachment of the same index.
Thus, the cells of $K_X$ correspond bijectively to the handles of the same index, and the claim follows from \cite[Corollary~8.3.14]{BHS11}.

\end{proof}

\begin{Remark}
\label{rem:identities-among-relations}
When the handle decomposition has a single $0$-handle, $\Pi X$ is reduced and gives the group presentation
\[
\pi_1(X)\cong(\Pi X)_1/\Ima(\partial_2)
\cong\langle g_a\mid\partial_2(r_b)\rangle.
\]
Since each $\partial_3(s_c)$ satisfies $\partial_2\partial_3(s_c)=1$, it is an identity among relations for this presentation \cite{Pei49,Rei49,BH82,Sie86}.
An element of $\Ker(\partial_2)$ records an identity sequence up to Peiffer equivalence.
If $K_X$ has several $0$-cells, choose a basepoint and a maximal tree in $K_X^{(1)}$, and collapse the tree to obtain a reduced model.
The resulting degree-$1$ generators and degree-$2$ boundaries give a presentation of the same form, and the images of the elements $\partial_3(s_c)$ in the reduced model give identities among its relations.
Every knot exterior admits a handle decomposition with one $0$-handle and no $3$-handles.
If this decomposition has $g$ $1$-handles, then $\chi(E_K)=0$ implies that it has $g-1$ $2$-handles.
The relative handle decomposition corresponding to the filling solid torus in Proposition~\ref{prop:crossed-complex-of-the-filling-solid-torus} adds one $2$-handle and one $3$-handle, and hence produces a Heegaard decomposition of genus $g$.
In the geometric presentation associated to a Heegaard diagram, Sieradski's construction represents the attaching sphere by a taut identity sequence \cite[pp.~132--133]{Sie86}.

\end{Remark}

\section{Boundary Torus Classes and the Dehn Filling Formula}
\label{sec:boundary-torus-classes-for-a-distinguished-boundary-component}

This section determines the generators and their boundaries that are added to the fundamental crossed complex when a distinguished boundary torus of a link exterior is Dehn filled along a slope $p/q$.
Definition~\ref{def:boundary-torus-class} records a cellular representative of the boundary component as a boundary torus class.
Lemma~\ref{lem:torus-conjugator-for-filling} constructs a slope word $X_{p/q}$, a complementary word $Y_{p/q}$, and a conjugating element $G_{p/q}$ from a continued fraction expansion and identifies their respective roles.
Proposition~\ref{prop:crossed-complex-of-the-filling-solid-torus} uses these three words to give a CW model of the filling solid torus relative to its boundary.
Finally, Theorem~\ref{thm:uniform-substitution-formula-for-p-over-q-surgery} glues this relative model to the link exterior and gives the Dehn filling formula.

We first fix a link exterior and a distinguished boundary component.
Let $L\subset S^3$ be an oriented link, and let $\nu L$ be its closed tubular neighborhood.
Define the link exterior by
\[
E_L:=S^3\setminus\operatorname{int}(\nu L).
\]
Then $E_L$ is a compact connected oriented $3$-manifold with torus boundary.
Choose a component $L_*\subset L$, and write $T_*:=\partial\nu(L_*)\subset\partial E_L$ for the corresponding boundary torus.
For the computation, fix a filtered CW model $K_L$ associated with a finite handle decomposition of $E_L$ and a homotopy equivalence $h_L\colon K_L\to E_L$, and write $\Pi E_L:=\Pi K_L$ for its fundamental crossed complex.

\subsection{Definition of Boundary Torus Classes}
\label{subsec:the-standard-torus-model}

To represent boundary parametrizations by cellular maps, we fix a product CW structure on a standard torus.
Give each factor circle the CW structure with $0$-cell $1$ and one oriented open $1$-cell.
The product CW structure on the oriented torus $T=S^1_\mu\times S^1_\lambda$ then has the unique $0$-cell $x_T=(1,1)$, the oriented open $1$-cells $e_\mu=(S^1_\mu\setminus\{1\})\times\{1\}$ and $e_\lambda=\{1\}\times(S^1_\lambda\setminus\{1\})$, and the oriented open $2$-cell $e_\Theta=(S^1_\mu\setminus\{1\})\times(S^1_\lambda\setminus\{1\})$.
Write $\mu,\lambda$ for the degree-$1$ generators of $\Pi T$ determined by the characteristic maps of $e_\mu,e_\lambda$, and write $\Theta$ for the degree-$2$ generator determined by $e_\Theta$.
We use the commutator convention $[a,b]:=aba^{-1}b^{-1}$.
With this CW structure, the boundary of $\Theta$ is
\[
\partial_2(\Theta)=[\mu,\lambda].
\]

We next fix a boundary parametrization and the framing conventions that determine it.
Choose a basepoint $b_*$ on $T_*$.
Orient a normal disk of $\nu(L_*)$ so that its orientation followed by the orientation of $L_*$ agrees with the orientation of $S^3$, and denote the boundary of this disk by $m_*$.
Fix a framing of $L_*$, orient its push-off in the same direction as $L_*$, and denote this push-off by $l_*$.
Here a framing means a homotopy class of trivializations of the oriented normal bundle of $L_*$.
Let $l_*^0$ be the preferred longitude with linking number $0$ with $L_*$.
There is then a unique $\kappa\in\Z$ such that $[l_*]=\kappa[m_*]+[l_*^0]$ in $H_1(T_*;\Z)$, and $\kappa=0$ is allowed.
Thus, we impose no nonzero condition on the integer framing coefficient.
Fix a basepoint-preserving homeomorphism
\[
s_0\colon(T,x_T)\xrightarrow{\cong}(T_*,b_*)
\]
that maps the factor circles $S^1_\mu$ and $S^1_\lambda$ of $T$ to $m_*$ and $l_*$, respectively.
When $T_*$ is oriented as the boundary of $E_L$, the ordered pair $(m_*,l_*)$ induces the opposite orientation, so $\deg s_0=-1$.
Fix the boundary parametrization $s$ as the composite
\[
s\colon T\xrightarrow{s_0}T_*\hookrightarrow E_L.
\]
If the chosen framing satisfies $l_*=l_*^0$, we call $s$ a preferred parametrization.
We measure all slopes below with respect to the ordered basis $(m_*,l_*)$.
Under the slope convention of Subsection~\ref{subsec:p-over-q-surgery}, the slope $1/0$ represents the meridian $m_*$.
It is therefore allowed as a Dehn filling slope but is not a longitude arising from a framing.

\begin{Definition}
\label{def:boundary-torus-class}
A boundary torus class with respect to the fixed boundary parametrization $s$ is a triple $(\Omega,\meri,\longi)$ for which there exist a filtration-preserving cellular map
\[
f_\Omega\colon T\longrightarrow K_L
\]
and a homotopy $H_\Omega\colon h_L\circ f_\Omega\simeq s$ such that
\[
\meri=(f_\Omega)_*(\mu),
\qquad
\longi=(f_\Omega)_*(\lambda),
\qquad
\Omega=(f_\Omega)_*(\Theta).
\]
If $x_0:=f_\Omega(x_T)$, then
\[
\meri,\longi\in(\Pi E_L)_1(x_0),
\qquad
\Omega\in(\Pi E_L)_2(x_0).
\]
The path $\gamma_\Omega(t):=H_\Omega(x_T,t)$ joins the image $h_L(x_0)$ of the basepoint of $K_L$ to the boundary basepoint $b_*$ and serves as the basepoint path used to compare them.

\end{Definition}

Boundary torus classes always exist and satisfy a commutator boundary formula.
Choose a homotopy inverse $k_L$ of $h_L$.
By cellular approximation, choose a filtration-preserving cellular map $f_\Omega$ homotopic to $k_L\circ s$ \cite{Bro06}.
Then $h_L\circ f_\Omega\simeq h_L\circ k_L\circ s\simeq s$, which proves existence.
Naturality of the boundary map also gives
\[
\partial_2(\Omega)=\partial_2((f_\Omega)_*(\Theta))=(f_\Omega)_*(\partial_2\Theta)=[\meri,\longi].
\]

\begin{Remark}
\label{rem:meridional-dehn-twist}
If the framing is shifted by an integer $k$, the longitude $l_*$ is replaced by $m_*^k l_*$, and the boundary torus class transforms as follows.
Let $d_k\colon T\to T$ be the based Dehn twist about the meridian that represents this change on the standard torus.
Define $D_k$ on the $1$-skeleton of $T$ by the edge words
\[
(D_k)_*(\mu)=\mu,
\qquad
(D_k)_*(\lambda)=\mu^k\lambda.
\]
Their commutator is null-homotopic in $T$, so this map extends over the $2$-cell to a cellular map $D_k\colon T\to T$.
Since $T$ is aspherical and $D_k$ induces the same map on $\pi_1(T)$ as $d_k$, the two maps are based homotopic.
The vanishing of $\pi_2(T)$ also makes $\partial_2\colon(\Pi T)_2\to(\Pi T)_1$ injective.
Moreover,
\[
\partial_2\bigl((D_k)_*(\Theta)\bigr)
=[\mu,\mu^k\lambda]
=\partial_2(\mu^k\cdot\Theta),
\]
and hence $(D_k)_*(\Theta)=\mu^k\cdot\Theta$.
Therefore, if $f_\Omega$ represents the boundary torus class $(\Omega,\meri,\longi)$ with respect to the boundary parametrization $s$, then the triple determined by $f_\Omega\circ D_k$ is
\[
(\meri^k\cdot\Omega,\ \meri,\ \meri^k\longi).
\]
This composite represents the boundary parametrization $s\circ d_k$ obtained by changing the framing by $k$.

\end{Remark}

\subsection{Slope Words from Continued Fractions}
\label{subsec:p-over-q-surgery}

\begin{Definition}
\label{def:continued-fraction-slope-data}
To represent a filling curve by a word in a free group, we now fix conventions for slopes and continued fraction expansions.
A slope is an unoriented primitive class in $H_1(T;\Z)\cong\Z^2$.
With respect to the ordered basis $([\mu],[\lambda])$ of $H_1(T;\Z)\cong\Z^2$, write $p/q$ for the finite slope represented by $p[\mu]+q[\lambda]$, where $\gcd(p,q)=1$ and $q>0$.
Write the infinite slope as $1/0$ and represent it by the column vector $(1,0)^T$ corresponding to $[\mu]$.
Thus, every slope has a unique expression $p/q$ such that $\gcd(p,q)=1$, $q\geq0$, and $p=1$ if $q=0$.
Define regular continued fractions by
\[
[c_1]=c_1,
\qquad
[c_1,\ldots,c_n]
=c_1+\frac{1}{[c_2,\ldots,c_n]}.
\]
For a finite slope $p/q$, take the unique regular continued fraction expansion
\[
p/q=[c_1,\ldots,c_n]
\]
such that $c_1\in\Z$, $c_i\in\Z_{>0}$ for $i\geq2$, and $c_n\geq2$ when $n>1$ \cite[Chapter~X]{HW08}.
If $n$ is odd, replace the last term $c_n$ by the two terms $c_n-1,1$.
This procedure uniquely determines an expansion of even length
\[
p/q=[a_1,\ldots,a_{2N}].
\]
For every finite slope $p/q$, we use this uniquely determined even-length expansion below.
The words associated with the infinite slope $1/0$ will be defined separately, without using the matrix recursion.

The matrix recursion used in the continued fraction construction is defined as follows.
Set
\[
R_a=
\begin{pmatrix}1&a\\0&1\end{pmatrix},
\qquad
L_a=
\begin{pmatrix}1&0\\a&1\end{pmatrix}.
\]
Let $M_0=I$, and for $1\leq i\leq 2N$, define
\[
M_i:=
\begin{cases}
M_{i-1}R_{a_i},&i\text{ is odd},\\
M_{i-1}L_{a_i},&i\text{ is even}.
\end{cases}
\]
Write $A_i$ and $B_i$ for the first and second columns of $M_i$, respectively, so that $M_i=(A_i\mid B_i)$.
See \cite[Chapter~X]{HW08} for the matrix computation of the convergents.

The next construction lifts the recursion defining the matrices $M_i$ to words in the free group carried by the $1$-skeleton of the standard torus; for Christoffel bases of a rank-two free group, see \cite[Section~3, especially Theorem~3.2]{KR07}.
The $1$-skeleton $T^{(1)}$ is the wedge of the $1$-cells $e_\mu$ and $e_\lambda$.
Therefore, the vertex group in degree $1$ of $\Pi T$ is
\[
(\Pi T)_1(x_T)=\pi_1(T^{(1)},x_T)=F(\mu,\lambda).
\]
For a word $W\in F(\mu,\lambda)$, let $[W]_{\mathrm{ab}}\in\Z^2$ be the column vector of exponent sums of $W$ with respect to the ordered basis $([\mu],[\lambda])$.
We call a free basis $(X,Y)$ of $F(\mu,\lambda)$ positively oriented if the matrix with columns $[X]_{\mathrm{ab}}$ and $[Y]_{\mathrm{ab}}$ has determinant $1$.
In $F(\mu,\lambda)$, set
\[
(X_0,Y_0,G_0):=(\mu,\lambda,1).
\]
For $i=1,\ldots,2N$, define
\[
(X_i,Y_i,G_i):=
\begin{cases}
\bigl(X_{i-1},\,X_{i-1}^{a_i}Y_{i-1},\,
X_{i-1}^{a_i}G_{i-1}\bigr),
&i\ \text{odd},\\[2mm]
\bigl(X_{i-1}Y_{i-1}^{a_i},\,Y_{i-1},\,G_{i-1}\bigr),
&i\ \text{even}.
\end{cases}
\]
Define the three final words associated with a finite slope $p/q$ by
\[
X_{p/q}:=X_{2N},
\qquad
Y_{p/q}:=Y_{2N},
\qquad
G_{p/q}:=G_{2N}.
\]
We call $X_{p/q}$ the signed Christoffel word of slope $p/q$.
For negative slopes its based representative may differ by conjugation from the convention in \cite[(3.1)]{KR07}.
The word $Y_{p/q}$ represents a complementary slope, and $G_{p/q}$ is a conjugating element relating $[X_{p/q},Y_{p/q}]$ to the standard commutator $[\mu,\lambda]$.
For the slope $1/0$, set
\[
X_{1/0}:=\mu,
\qquad
Y_{1/0}:=\lambda,
\qquad
G_{1/0}:=1.
\]

\end{Definition}

\begin{Lemma}
\label{lem:torus-conjugator-for-filling}
For every slope $p/q$, including $1/0$, the words $X_{p/q}$, $Y_{p/q}$, and $G_{p/q}$ satisfy the following properties.
\begin{enumerate}[(1)]
\item $(X_{p/q},Y_{p/q})$ is a positively oriented free basis of $F(\mu,\lambda)$.
\item $\bigl[X_{p/q}\bigr]_{\mathrm{ab}}=(p,q)^T$.
\item $\partial_2(G_{p/q}\cdot\Theta)=[X_{p/q},Y_{p/q}]$.
\end{enumerate}

\end{Lemma}

\begin{proof}
We first prove the three claims for a finite slope $p/q=[a_1,\ldots,a_{2N}]$ by induction over the stages of the recursion.
Multiplying $M_{i-1}=(A_{i-1}\mid B_{i-1})$ on the right by $R_{a_i}$ or $L_{a_i}$ changes its two columns according to
\[
\bigl(A_i\mid B_i\bigr)=
\begin{cases}
\bigl(A_{i-1}\mid B_{i-1}+a_iA_{i-1}\bigr),&i\text{ is odd},\\
\bigl(A_{i-1}+a_iB_{i-1}\mid B_{i-1}\bigr),&i\text{ is even}.
\end{cases}
\]
The recursion for $(X_{i-1},Y_{i-1})$ induces the same updates of the two columns after abelianization.
Induction from the initial values therefore gives
\[
\bigl[X_i\bigr]_{\mathrm{ab}}=A_i,
\qquad
\bigl[Y_i\bigr]_{\mathrm{ab}}=B_i.
\]
Each update $(X,Y)\mapsto(X,X^aY)$ or $(X,Y)\mapsto(XY^a,Y)$ is an invertible Nielsen transformation, and the induced matrix on the abelianization is $R_a$ or $L_a$, respectively.
Both matrices have determinant $1$, so every pair $(X_i,Y_i)$ is a positively oriented free basis.
Using the identities $[X,X^aY]=X^a[X,Y]X^{-a}$ and $[XY^a,Y]=[X,Y]$ together with the recursion for $G_i$, the same induction gives
\[
[X_i,Y_i]=G_i[\mu,\lambda]G_i^{-1}.
\]

The final matrix identifies the abelianized slope word and completes the boundary computation.
Set $C_a:=\left(\begin{smallmatrix}a&1\\1&0\end{smallmatrix}\right)$.
Since $R_aL_b=C_aC_b$, we have
\[
M_{2N}
=R_{a_1}L_{a_2}\cdots R_{a_{2N-1}}L_{a_{2N}}
=C_{a_1}\cdots C_{a_{2N}}.
\]
The recursion for the convergents of a regular continued fraction \cite[Chapter~X]{HW08} shows that the first column of this product is $(p,q)^T$.
Hence $\bigl[X_{p/q}\bigr]_{\mathrm{ab}}=(p,q)^T$.
Compatibility of the crossed module action with the boundary map gives
\[
\partial_2(G_{p/q}\cdot\Theta)
=G_{p/q}\partial_2(\Theta)G_{p/q}^{-1}
=G_{p/q}[\mu,\lambda]G_{p/q}^{-1}
=[X_{p/q},Y_{p/q}].
\]
For the slope $1/0$, all three conclusions follow immediately from $(X_{1/0},Y_{1/0},G_{1/0})=(\mu,\lambda,1)$.

\end{proof}

\subsection{The Crossed Complex of the Filling Solid Torus}

We describe the filling solid torus as a CW complex relative to the boundary torus $T$.
First, we attach a $2$-cell to $T$ to kill the loop representing the slope $p/q$ and compute the second homotopy group of the resulting $2$-skeleton $K_2$.
We then attach a $3$-cell along a free generator of this group to construct a model $V_{p/q}$ that is homotopy equivalent to the actual filling solid torus relative to $T$.

We construct a $2$-skeleton adapted to the slope.
For brevity, write the three words constructed in the preceding subsection as
\[
x:=X_{p/q},
\qquad
y:=Y_{p/q},
\qquad
\Theta_0:=G_{p/q}\cdot\Theta.
\]
Lemma~\ref{lem:torus-conjugator-for-filling} gives
\[
\partial_2(\Theta_0)=[x,y]
\]
for every slope $p/q$.
Attach a $2$-cell $e_r^2$ to $T$ along the based cellular loop represented by $x$, and set
\[
K_2:=T\cup_x e^2_r.
\]
If $r$ denotes the generator of $(\Pi K_2)_2$ represented by this $2$-cell, then $\partial_2(r)=x$.

We next compare $K_2$ with a model whose $2$-cell is attached along the embedded slope.
Let
\[
c_{p/q}\colon(S^1,1)\hookrightarrow(T,x_T)
\]
be the oriented embedded slope represented by $(p,q)$, and set $K_2^c:=T\cup_{c_{p/q}}e^2$.
This comparison replaces the word $x$ obtained from the continued fraction with the geometric slope $c_{p/q}$ along which the meridian disk of the solid torus is attached.
Lemma~\ref{lem:torus-conjugator-for-filling} gives $[x]_{\mathrm{ab}}=(p,q)^T$.
The loop $c_{p/q}$ represents the same element of $\pi_1(T,x_T)\cong\Z^2$, so the two based loops are homotopic in $T$ relative to the basepoint.
Homotopy invariance of attaching maps therefore allows us to choose a homotopy equivalence relative to $T$,
\[
F\colon K_2\xrightarrow{\simeq}K_2^c,
\qquad
F|_T=\mathrm{id}_T.
\]

We define a candidate for the boundary of the $3$-cell attached to $K_2$ by
\begin{equation}
\label{eq:eta-solid-torus-generator}
\eta
:=
r^{-1}\Theta_0(y\cdot r)
=
r^{-1}(G_{p/q}\cdot\Theta)(Y_{p/q}\cdot r).
\end{equation}
The commutator convention gives
\[
\partial_2(\eta)=x^{-1}[x,y](yxy^{-1})=1.
\]
Thus,
\[
\eta\in\Ker\bigl(\partial_2\colon(\Pi K_2)_2\longrightarrow(\Pi K_2)_1\bigr),
\]
and $\eta$ represents a spherical class in $K_2$.

\begin{Lemma}
\label{lem:pi-two-after-slope-two-handle}
The second homotopy group of $K_2$ is freely generated by $\eta$ over its integral group ring.
More precisely, let $R:=\Z[\pi_1(K_2)]$.
Then
\[
\pi_1(K_2)\cong\langle y\rangle\cong\Z,
\qquad
\pi_2(K_2)\cong R\eta,
\]
and $\eta$ is a free $R$-generator.

\end{Lemma}

\begin{proof}
We first identify the fundamental group and the homotopy class represented by $\eta$.
By Lemma~\ref{lem:torus-conjugator-for-filling}, $(x,y)$ is a free basis, and hence
\[
\pi_1(K_2)=\langle x,y\mid[x,y],x\rangle
\cong\langle y\rangle.
\]
As computed immediately after \eqref{eq:eta-solid-torus-generator}, $\partial_2(\eta)=1$.
Since $K_2^{(1)}$ is a graph, $\pi_2(K_2^{(1)})=0$.
The relative homotopy exact sequence of $(K_2,K_2^{(1)})$ therefore identifies $\pi_2(K_2)$ with the kernel of $\partial_2$.
Thus, $\eta$ represents an element of $\pi_2(K_2)$.

We compute $\pi_2(K_2)$ from the cellular chain complex of the universal cover.
Let \(\widetilde K_2\to K_2\) be the universal cover, and regard its cellular chain groups as left modules over \(R=\Z[\pi_1(K_2)]\) through the action of the deck transformation group \(\pi_1(K_2)=\langle y\rangle\).
Write $t$ for the unit of the group ring $R$ corresponding to the generator $y\in\pi_1(K_2)$.
As observed in the proof of Lemma~\ref{lem:torus-conjugator-for-filling}, a sequence of orientation-preserving Nielsen transformations carries $(\mu,\lambda)$ to $(x,y)$.
Fix a lift of the basepoint $x_T$, and denote by $e_x$ and $e_y$ the cellular $1$-chains obtained by lifting the loops $x$ and $y$ from this point.
Each elementary Nielsen transformation induces an invertible $R$-linear change of basis on the lifted $1$-chains, and the inverse Nielsen transformation induces its inverse.
Therefore, $(e_x,e_y)$ is an $R$-basis of $C_1(\widetilde K_2)$.
Moreover, $\Theta_0=G_{p/q}\cdot\Theta$ is a unit multiple of a lift of the boundary torus $2$-cell $\Theta$, so $(\Theta_0,r)$ is an $R$-basis of $C_2(\widetilde K_2)$.

As in Trotter's calculation based on group presentations, the cellular boundary of a lifted $2$-cell is determined by the Fox derivatives of its attaching word \cite{Tro62}.
Let
\[
\phi\colon\Z[F(x,y)]\longrightarrow R,
\qquad
\phi(x)=1,
\quad
\phi(y)=t
\]
be the homomorphism of group rings induced by the quotient map $F(x,y)\to\pi_1(K_2)$.
Under our convention of left $R$-modules,
\[
\phi\left(\frac{\partial[x,y]}{\partial x}\right)=1-t,
\qquad
\phi\left(\frac{\partial[x,y]}{\partial y}\right)=0,
\qquad
\phi\left(\frac{\partial x}{\partial x}\right)=1,
\qquad
\phi\left(\frac{\partial x}{\partial y}\right)=0.
\]
With respect to the bases above, the chain groups and the boundary map are
\nopagebreak[4]
\[
\begin{gathered}
C_2(\widetilde K_2)=R\Theta_0\oplus Rr,
\qquad
C_1(\widetilde K_2)=Re_x\oplus Re_y,\\
d_2(\Theta_0)=(1-t)e_x,
\quad
d_2(r)=e_x.
\end{gathered}
\]
Since $K_2$ is $2$-dimensional, $H_2(\widetilde K_2)=\Ker d_2$.
The boundary formulas show that this kernel is the free $R$-module generated by $\Theta_0+(t-1)r$.
Hence,
\[
H_2(\widetilde K_2)=R\{\Theta_0+(t-1)r\}.
\]

We now identify $\eta$ with the generator obtained from the calculation in the universal cover.
The covering map induces an isomorphism $\pi_2(\widetilde K_2)\to\pi_2(K_2)$, and the Hurewicz isomorphism $\pi_2(\widetilde K_2)\cong H_2(\widetilde K_2)$ holds because $\widetilde K_2$ is simply connected.
Denote the composite of the inverse of the first isomorphism with the Hurewicz isomorphism by
\[
\mathfrak h\colon\pi_2(K_2)\xrightarrow{\cong}H_2(\widetilde K_2).
\]
Replacing the product in \eqref{eq:eta-solid-torus-generator} by a sum of cellular $2$-chains and the action of $y$ by the deck transformation $t$, we find that $\mathfrak h(\eta)$ is represented by the cellular $2$-cycle
\[
-r+\Theta_0+tr
=\Theta_0+(t-1)r.
\]
This cycle is the free $R$-generator of $H_2(\widetilde K_2)$ found above.
Therefore, $\eta$ is a free $R$-generator of $\pi_2(K_2)$.

\end{proof}

We next relate the algebraic model to the geometric filling solid torus.
Let $\widehat V_{p/q}$ be a solid torus, and choose a boundary homeomorphism
\[
\widehat j_{p/q}\colon T\xrightarrow{\cong}\partial\widehat V_{p/q}
\]
such that $\widehat j_{p/q}\circ c_{p/q}$ bounds an oriented meridian disk.
We use $\widehat j_{p/q}$ to regard $T$ as a boundary subspace of $\widehat V_{p/q}$.
Recall that $K_2^c=T\cup_{c_{p/q}}e^2$.
The standard relative handle decomposition of a solid torus, together with the standard correspondence from handle decompositions to CW decompositions, shows that $\widehat V_{p/q}$ is homotopy equivalent relative to $T$ to a CW complex obtained by attaching one $3$-cell to $K_2^c$.

\begin{Proposition}
\label{prop:crossed-complex-of-the-filling-solid-torus}
The filling solid torus admits the following filtered CW model relative to its boundary torus.
There exists a filtered CW pair $(V_{p/q},T)$ equipped with the inclusion of the boundary subcomplex $j_{p/q}\colon T\hookrightarrow V_{p/q}$ such that $(V_{p/q},T)$ is homotopy equivalent to $(\widehat V_{p/q},T)$ relative to $T$.
The crossed complex $\Pi V_{p/q}$ is obtained from $\Pi T$ by adjoining a generator $r$ in degree $2$ and a generator $\sigma$ in degree $3$.
The whisker and orientation of $\sigma$ can be chosen so that
\begin{equation}
\label{eq:solid-torus-three-boundary-formula}
\partial_2(r)=X_{p/q},
\qquad
\partial_3(\sigma)
=r^{-1}(G_{p/q}\cdot\Theta)(Y_{p/q}\cdot r).
\end{equation}

\end{Proposition}

The following diagram summarizes the relation between the two CW complexes compared in the proposition.
The hooked horizontal arrows are cellular inclusions obtained by attaching the indicated $3$-cells.
The vertical arrows and the final horizontal arrow are homotopy equivalences relative to $T$.
\begin{center}
\small
\renewcommand{\arraystretch}{1.35}
\(
\begin{array}{c@{\;}c@{\;}c@{\;}c}
K_2=T\cup_x e_r^2
&
\xhookrightarrow{\quad\text{attach a }3\text{-cell}\quad}
&
V_{p/q}
&
\\[3pt]
\Big\downarrow\!\mathrlap{\scriptstyle\,F\;(\simeq_T)}
&&
\Big\downarrow\!\mathrlap{\scriptstyle\,\simeq_T}
&\\[3pt]
K_2^c=T\cup_{c_{p/q}}e^2
&
\xhookrightarrow{\quad\text{attach the standard }3\text{-cell}\quad}
&
K_2^c\cup e^3
&
\xrightarrow{\;\simeq_T\;}\widehat V_{p/q}
\end{array}
\)
\end{center}
The diagram is not asserted to commute strictly.
It is homotopy commutative in the sense that the image under $F$ of the upper attaching map is homotopic to the lower attaching map.
In the crossed complexes, denote by $r$ and $\sigma$ the generators represented by the $2$-cell in $K_2$ and the subsequently attached $3$-cell, respectively.
Their boundaries are $\partial_2(r)=x$ and $\partial_3(\sigma)=\eta$.

\begin{proof}
We first compare the actual attaching class with the model above.
Let $z^c\colon S^2\to K_2^c$ be the attaching map of the $3$-cell in the standard relative handle decomposition, and let $z\colon S^2\to K_2$ represent $F_*^{-1}([z^c])\in\pi_2(K_2)$.
Then $F\circ z\simeq z^c$.
By homotopy invariance of attaching maps \cite[Corollary~1 to Theorem~7.5.5]{Bro06} and the gluing theorem for adjunction spaces \cite[Theorem~7.5.7]{Bro06}, the map $F$ extends to a homotopy equivalence relative to $T$,
\[
\overline F\colon
K_2\cup_z e^3
\xrightarrow{\simeq_T}
K_2^c\cup_{z^c}e^3.
\]
On the other hand, the correspondence between the standard relative handle decomposition of the solid torus and its CW decomposition gives
\[
K_2^c\cup_{z^c}e^3\simeq_T\widehat V_{p/q}.
\]
Therefore, $V_{p/q}:=K_2\cup_z e^3$ is homotopy equivalent to $\widehat V_{p/q}$ relative to $T$.
Write $j_{p/q}\colon T\hookrightarrow V_{p/q}$ for the inclusion of the boundary subcomplex in this model.

The attaching class $[z]$ generates $\pi_2(K_2)$ over the group ring.
Attaching one $3$-cell does not change the fundamental group, so $K_2\to V_{p/q}$ induces an isomorphism on $\pi_1$.
Moreover, because $V_{p/q}$ is obtained from $K_2$ by attaching only a $3$-cell, the pair $(V_{p/q},K_2)$ and the corresponding pair of universal covers are $2$-connected.
As in Lemma~\ref{lem:pi-two-after-slope-two-handle}, set $R:=\Z[\pi_1(K_2)]\cong\Z[t,t^{-1}]$.
The isomorphism $\pi_1(K_2)\cong\pi_1(V_{p/q})$ also identifies this ring with $\Z[\pi_1(V_{p/q})]$.
The pair of universal covers has one $R$-orbit of relative $3$-cells.
The equivariant relative Hurewicz theorem gives
\[
\pi_3(V_{p/q},K_2)
\cong H_3(\widetilde V_{p/q},\widetilde K_2)
\cong R,
\]
and the boundary of this relative generator is $[z]$.
Since $V_{p/q}\simeq\widehat V_{p/q}\simeq S^1$, the homotopy exact sequence of the pair gives
\[
\pi_2(K_2)=R[z].
\]

It remains to normalize the unit relating the two generators.
Lemma~\ref{lem:pi-two-after-slope-two-handle} gives $\pi_2(K_2)=R\eta$, and the equality $\pi_2(K_2)=R[z]$ therefore gives $R[z]=R\eta$.
Thus, if $[z]=a(t)\eta$, then $a(t)$ is a unit in $R$.
The only units of $R=\Z[t,t^{-1}]$ are $\pm t^n$, so
\[
[z]=\pm t^n\eta
\qquad(n\in\Z).
\]
The factor $t^n$ corresponds to changing the whisker of the $3$-cell, and the sign corresponds to reversing its orientation.
After changing the whisker and orientation of $\sigma$, we may therefore arrange that
\[
\partial_3(\sigma)=\eta
=r^{-1}(G_{p/q}\cdot\Theta)(Y_{p/q}\cdot r).
\]
The equality $\partial_2(\eta)=1$ shows that this element lies in the kernel of $\partial_2$.

\end{proof}

\begin{Remark}
\label{rem:filling-complement-independence}
We now compare the boundary elements obtained from different positively oriented complements.
If the complement is changed by an integral multiple of the slope vector and the associated word and conjugator are changed accordingly, the boundary element is unchanged.
Set $A:=[x]_{\mathrm{ab}}$ and $B:=[y]_{\mathrm{ab}}$.
By Lemma~\ref{lem:torus-conjugator-for-filling}, $\det(A,B)=1$, and every positively oriented complement $B'$ of $A$ satisfying $\det(A,B')=1$ has the unique form $B+kA$.
Set $B_k:=B+kA$.
A word with $[y_k]_{\mathrm{ab}}=B_k$ and its corresponding conjugator may be chosen as
\[
y_k:=x^ky,
\qquad
G_k:=x^kG_{p/q}.
\]
Using $\partial_2(r)=x$ and the Peiffer identity, we obtain
\[
\begin{aligned}
r^{-1}(G_k\cdot\Theta)(y_k\cdot r)
&=r^{-1}(x^k\cdot\Theta_0)\bigl((x^ky)\cdot r\bigr)\\
&=r^k\eta r^{-k}\\
&=\eta.
\end{aligned}
\]
The second equality uses the Peiffer identity $(\partial_2r)^k\cdot w=r^kwr^{-k}$, and the third uses $\eta\in\Ker(\partial_2)\subseteq Z((\Pi K_2)_2)$.

\end{Remark}

\subsection{The Substitution Formula for Dehn Filling}

We now construct a filtered CW model of the actual Dehn filling and describe its crossed complex explicitly from the boundary torus class.
Define the actual $p/q$-Dehn filling by
\[
E_L(p/q):=E_L\cup_{s,\,\widehat j_{p/q}}\widehat V_{p/q}.
\]
Here the gluing symbol denotes the pushout obtained from the disjoint union $E_L\sqcup\widehat V_{p/q}$ by identifying the two images $s(a)$ and $\widehat j_{p/q}(a)$ for every $a\in T$; explicitly, it is
\[
(E_L\sqcup\widehat V_{p/q})/
\{s(a)\sim\widehat j_{p/q}(a)\mid a\in T\}.
\]
Fix a boundary torus class $(\Omega,\meri,\longi)$ and write the slope in the normalized form $p/q$ fixed above, allowing $1/0$.
Also fix a representative $f_\Omega\colon T\to K_L$ of this boundary torus class and a comparison homotopy $H_\Omega\colon h_L\circ f_\Omega\simeq s$.
Write
\[
j_{p/q}\colon T\hookrightarrow V_{p/q}
\]
for the inclusion of the boundary subcomplex obtained in Proposition~\ref{prop:crossed-complex-of-the-filling-solid-torus}.
Using this inclusion and $f_\Omega$ as the gluing maps, define
\[
\mathcal X_{p/q}(f_\Omega)
:=
K_L\cup_{f_\Omega,j_{p/q}}V_{p/q}.
\]

The actual Dehn filling and its CW model are represented by the following pushout diagrams.
\begin{center}
\small
\renewcommand{\arraystretch}{1.35}
\(
\begin{array}{c@{\qquad\qquad}c}
\begin{array}{ccc}
T
& \xrightarrow{\;j_{p/q}\;}
& V_{p/q}
\\[3pt]
{\scriptstyle f_\Omega}\!\Big\downarrow
&&
\Big\downarrow
\\[3pt]
K_L
& \longrightarrow
& \mathcal X_{p/q}(f_\Omega)
\end{array}
&
\begin{array}{ccc}
T
& \xrightarrow{\;\widehat j_{p/q}\;}
& \widehat V_{p/q}
\\[3pt]
{\scriptstyle s}\!\Big\downarrow
&&
\Big\downarrow
\\[3pt]
E_L
& \longrightarrow
& E_L(p/q)
\end{array}
\\[5pt]
\text{Pushout for the CW model}
&
\text{Pushout for the actual filling}
\end{array}
\)
\end{center}
The left-hand square represents gluing by $f_\Omega$ and $j_{p/q}$, and the right-hand square represents gluing by $s$ and $\widehat j_{p/q}$.

We compare the CW model with the actual filling through three homotopy equivalences.
The map $h_L\colon K_L\to E_L$ is a homotopy equivalence, and $j_{p/q}$ is a closed cofibration.
The gluing theorem for adjunction spaces \cite[Theorem~7.5.7]{Bro06} therefore gives the first homotopy equivalence below.
The comparison homotopy $H_\Omega\colon h_L\circ f_\Omega\simeq s$ and homotopy invariance of attaching maps \cite[Corollary~1 to Theorem~7.5.5]{Bro06} give the second.
Finally, gluing with the homotopy equivalence relative to $T$ from Proposition~\ref{prop:crossed-complex-of-the-filling-solid-torus} gives the third.
Thus,
\[
\begin{aligned}
\mathcal X_{p/q}(f_\Omega)
&\simeq E_L\cup_{h_L\circ f_\Omega,\,j_{p/q}}V_{p/q}\\
&\simeq E_L\cup_{s,\,j_{p/q}}V_{p/q}\\
&\simeq E_L\cup_{s,\,\widehat j_{p/q}}\widehat V_{p/q}
=E_L(p/q).
\end{aligned}
\]

We equip the gluing model with a filtration by defining its skeleton in each degree $n$ as
\[
\mathcal X_{p/q}(f_\Omega)^{(n)}
:=
K_L^{(n)}
\cup_{f_\Omega|_{T^{(n)}},\,j_{p/q}|_{T^{(n)}}}
V_{p/q}^{(n)}.
\]
The map $j_{p/q}|_{T^{(n)}}$ is the inclusion of a subcomplex and a closed cofibration, and $f_\Omega|_{T^{(n)}}$ is cellular.
This filtration therefore gives $\mathcal X_{p/q}(f_\Omega)$ a CW structure with a connected filtration.
By the pushout theorem for crossed complexes \cite[Theorem~8.2.5]{BHS11},
\[
\Pi\mathcal X_{p/q}(f_\Omega)
\cong
\Pi E_L\sqcup_{\Pi T}\Pi V_{p/q}
\]
is a pushout in the category $\Crs$ of crossed complexes.
Its structure maps are
\[
(f_\Omega)_*\colon\Pi T\longrightarrow\Pi E_L,
\qquad
(j_{p/q})_*\colon\Pi T\longrightarrow\Pi V_{p/q}.
\]

We introduce the substitution notation used in the uniform Dehn filling formula.
Let $W(\mu,\lambda)\in F(\mu,\lambda)$.
For loops $a,b\in C_1(z)$ based at the same object $z$ of a crossed complex $C$, let $W(a,b)$ denote the element obtained from the word $W$ by substituting $\mu\mapsto a$ and $\lambda\mapsto b$.

\begin{Theorem}
\label{thm:uniform-substitution-formula-for-p-over-q-surgery}
The uniform Dehn filling formula describes the crossed complex of the filled space as follows.
Fix a boundary torus class $(\Omega,\meri,\longi)$ and a slope $p/q$, and define elements of the vertex group $(\Pi E_L)_1(x_0)$ by
\[
\alpha:=X_{p/q}(\meri,\longi),
\qquad
\beta:=Y_{p/q}(\meri,\longi),
\qquad
G_\Omega:=G_{p/q}(\meri,\longi)=(f_\Omega)_*(G_{p/q}).
\]
Then $\Pi\mathcal X_{p/q}(f_\Omega)$ is the crossed complex obtained from $\Pi E_L$ by adjoining a generator $\rho$ in degree $2$ and a generator $\sigma$ in degree $3$, both based at the object $x_0$, and their boundaries satisfy
\begin{equation}
\label{eq:uniform-substitution-formula-for-p-over-q-surgery}
\partial_2(\rho)=\alpha,
\qquad
\partial_3(\sigma)
=\rho^{-1}(G_\Omega\cdot\Omega)(\beta\cdot\rho).
\end{equation}
\end{Theorem}

\begin{proof}
The formula follows by tracing the generators from the filling solid torus through the pushout.
For the map
\[
\Pi V_{p/q}\longrightarrow\Pi\mathcal X_{p/q}(f_\Omega),
\]
write $\rho$ for the image of the generator $r$ in Proposition~\ref{prop:crossed-complex-of-the-filling-solid-torus}, and retain the symbol $\sigma$ for the image of $\sigma$.
This map sends the generators as follows:
\[
\mu\mapsto\meri,
\qquad
\lambda\mapsto\longi,
\qquad
\Theta\mapsto\Omega,
\qquad
r\mapsto\rho.
\]
Consequently,
\[
X_{p/q}\mapsto\alpha,
\qquad
Y_{p/q}\mapsto\beta,
\qquad
G_{p/q}\cdot\Theta\mapsto G_\Omega\cdot\Omega.
\]
Taking the image of \eqref{eq:solid-torus-three-boundary-formula} in the pushout gives \eqref{eq:uniform-substitution-formula-for-p-over-q-surgery}.

\end{proof}

The same construction applies simultaneously to several boundary components.
For pairwise distinct boundary components $T_1,\ldots,T_n$, specify a boundary torus class $(\Omega_i,\meri_i,\longi_i)$ and a slope $p_i/q_i$ for each $i$.
Since the boundary components are disjoint, the relative models of the corresponding filling solid tori can be glued to $E_L$ simultaneously to form one filtered CW pushout.
Applying the pushout theorem for crossed complexes to this gluing adds one generator in degree $2$ and one in degree $3$ for each component, with their boundaries given by the same substitution formula as in Theorem~\ref{thm:uniform-substitution-formula-for-p-over-q-surgery}.

\section{Crossed Complexes of Satellite Exteriors}
\label{sec:torus-gluings-and-operations}

This section constructs a filtered CW model of a satellite exterior from boundary parametrizations of a pattern exterior $X_P$ and a companion exterior $E_K$.
We first define the satellite knot and the actual gluing map, and then check its orientation and the condition for extending it over the solid torus.
We then prove in Theorem~\ref{thm:satellite-pushout} that the crossed complex of the resulting model is computed by a pushout over the torus.

We first define the pattern and companion exteriors and fix their boundary coordinates.
Let $V$ be an oriented solid torus, let $P\subset\operatorname{int}(V)$ be a pattern knot, and set
\[
X_P:=V\setminus\operatorname{int}(\nu P).
\]
Write $\partial_{\mathrm{out}}X_P:=\partial V$ for the outer boundary and $\partial_{\mathrm{in}}X_P:=\partial\nu P$ for the inner boundary.
Let $K\subset S^3$ be a companion knot, and set
\[
E_K:=S^3\setminus\operatorname{int}(\nu K).
\]
Fix an orientation-preserving homeomorphism
\[
\Phi\colon V\xrightarrow{\cong}\nu K,
\]
and define the satellite knot to be the image $P(K):=\Phi(P)\subset S^3$.
Let
\[
s_{\mathrm{out}}\colon T\xrightarrow{\cong}\partial_{\mathrm{out}}X_P\hookrightarrow X_P
\]
be the parametrization of the outer boundary, and let
\[
s_K\colon T\xrightarrow{\cong}\partial E_K\hookrightarrow E_K
\]
be the boundary parametrization on the companion side.
When taking inverse maps and degrees, we omit the boundary inclusions and regard these parametrizations as homeomorphisms onto the corresponding boundary components.
On the domain $T$ of each parametrization, choose an ordered peripheral basis of $\pi_1(T)$ whose first element is the meridian.
Denote these bases by $(\mu_V,\lambda_V)$ and $(\mu_K,\lambda_K)$, and use the same notation for their images on the boundary components.

We now record the orientation and extension properties of the gluing map.
Identify $\partial\nu K$ and $\partial E_K$ as the same underlying space.
Their orientations are opposite because one is oriented as the boundary of $\nu K$ and the other as the boundary of $E_K$.
Transfer the actual gluing map $\Phi|_{\partial V}\colon\partial V\to\partial\nu K$ to the standard torus by defining
\[
\varphi_{\mathrm{top}}
:=s_K^{-1}\circ\Phi|_{\partial V}\circ s_{\mathrm{out}}
\colon T\longrightarrow T.
\]
The definition gives
$s_K\circ\varphi_{\mathrm{top}}=\Phi|_{\partial V}\circ s_{\mathrm{out}}$.
Since $\Phi$ preserves the orientations of $V$ and $\nu K$, the reversal of boundary orientations noted above implies that $\Phi|_{\partial V}$ reverses orientation as a gluing map from $\partial_{\mathrm{out}}X_P$ to $\partial E_K$.
Thus,
\[
\deg(\varphi_{\mathrm{top}})
=-\deg(s_{\mathrm{out}})\deg(s_K).
\]
Moreover, because $\varphi_{\mathrm{top}}$ is induced on the boundary by the homeomorphism $\Phi$ of solid tori, it satisfies
\[
(\varphi_{\mathrm{top}})_*[\mu_V]=\varepsilon[\mu_K]
\quad\text{in }H_1(T;\Z),
\qquad
\varepsilon\in\{\pm1\}.
\]
Conversely, a torus homeomorphism extends over the entire solid torus if and only if the image of the meridian has zero longitude component, or equivalently, if it sends the meridian to the meridian up to sign.
Since $\pi_1(T)$ is abelian, the same condition may be written as $\varphi_{\mathrm{top},*}(\mu_V)=\mu_K^\varepsilon$.

We construct a filtered gluing model of the satellite exterior.
The actual exterior of the satellite knot is the pushout obtained by identifying $\partial V$ with $\partial\nu K$ through $\Phi|_{\partial V}$:
\[
E_{P(K)}
=
X_P\cup_{s_{\mathrm{out}},\,s_K\circ\varphi_{\mathrm{top}}}E_K.
\]
Choose CW structures with connected filtrations on $X_P$ and $E_K$ such that, for every $n$,
\[
s_K(T^{(n)})=s_K(T)\cap E_K^{(n)},
\]
and such that $s_K|_{T^{(n)}}$ is the inclusion of a closed subcomplex.
Choose a filtration-preserving cellular map
\[
f_\varphi\colon T\longrightarrow X_P,
\qquad
f_\varphi\simeq s_{\mathrm{out}}\circ\varphi_{\mathrm{top}}^{-1}.
\]
Define
\[
Z_{\mathrm{sat}}^{(n)}
:=
X_P^{(n)}
\cup_{f_\varphi|_{T^{(n)}},\,s_K|_{T^{(n)}}}
E_K^{(n)}
\]
and set
\[
Z_{\mathrm{sat}}
:=
X_P\cup_{f_\varphi,\,s_K}E_K.
\]
The following diagram summarizes the boundary coordinates and the two gluings.
The copy of $T$ on the left is the standard torus with basis $(\mu_V,\lambda_V)$, and the copy on the right is the standard torus with basis $(\mu_K,\lambda_K)$.
The lower-left pushout for the actual exterior uses the standard torus on the left as its common subspace, whereas the lower-right pushout for the CW model uses the standard torus on the right.
\begin{center}
\small
\renewcommand{\arraystretch}{1.30}
\(
\begin{array}{c}
\begin{array}{ccc}
T
& \xrightarrow{\;\varphi_{\mathrm{top}}\;}
& T
\\[3pt]
{\scriptstyle s_{\mathrm{out}}}\!\Big\downarrow
&&
{\scriptstyle s_K}\!\Big\downarrow
\\[3pt]
\partial_{\mathrm{out}}X_P
& \xrightarrow{\;\Phi|_{\partial V}\;}
& \partial E_K
\end{array}
\\[5pt]
\text{Commutative diagram transferring the boundary coordinates to standard tori}
\\[9pt]
\begin{array}{c@{\qquad}c}
\begin{array}{ccc}
T
& \xrightarrow{\;s_K\circ\varphi_{\mathrm{top}}\;}
& E_K
\\[3pt]
{\scriptstyle s_{\mathrm{out}}}\!\Big\downarrow
&&
\Big\downarrow
\\[3pt]
X_P
& \longrightarrow
& E_{P(K)}
\end{array}
&
\begin{array}{ccc}
T
& \xrightarrow{\;s_K\;}
& E_K
\\[3pt]
{\scriptstyle f_\varphi}\!\Big\downarrow
&&
\Big\downarrow
\\[3pt]
X_P
& \longrightarrow
& Z_{\mathrm{sat}}
\end{array}
\\[5pt]
\text{Pushout for the actual exterior}
&
\text{Pushout for the CW model}
\\[4pt]
\multicolumn{2}{c}{
\displaystyle
f_\varphi\simeq
s_{\mathrm{out}}\circ\varphi_{\mathrm{top}}^{-1}
}
\end{array}
\end{array}
\)
\end{center}

\begin{Theorem}
\label{thm:satellite-pushout}
The satellite exterior satisfies the following torus-gluing formula.
The space $Z_{\mathrm{sat}}$ is a filtered CW complex, and
\[
\Pi Z_{\mathrm{sat}}
\cong
\Pi X_P\sqcup_{\Pi T}\Pi E_K
\]
holds in $\Crs$.
The two structure maps of this pushout are
\[
(f_\varphi)_*\colon\Pi T\longrightarrow\Pi X_P,
\qquad
(s_K)_*\colon\Pi T\longrightarrow\Pi E_K.
\]
Moreover, $Z_{\mathrm{sat}}\simeq E_{P(K)}$.

\end{Theorem}

\begin{proof}
We prove both the crossed-complex pushout formula and the identification with the actual exterior.
Since $s_K|_{T^{(n)}}$ is a closed cofibration and $f_\varphi|_{T^{(n)}}$ is cellular, every $Z_{\mathrm{sat}}^{(n)}$ is a CW complex.
The pushout theorem of Brown--Higgins \cite{BH81,BHS11} gives
\[
\Pi Z_{\mathrm{sat}}
\cong
\Pi X_P\sqcup_{\Pi T}\Pi E_K.
\]
Since $s_K$ is a cofibration and $f_\varphi\simeq s_{\mathrm{out}}\circ\varphi_{\mathrm{top}}^{-1}$,
\[
\begin{aligned}
Z_{\mathrm{sat}}
&\simeq
X_P\cup_{s_{\mathrm{out}}\circ\varphi_{\mathrm{top}}^{-1},\,s_K}E_K\\
&\cong
X_P\cup_{s_{\mathrm{out}},\,s_K\circ\varphi_{\mathrm{top}}}E_K\\
&=E_{P(K)}.
\end{aligned}
\]

\end{proof}

\section{Computing Boundary Torus Classes from Annulus Decompositions}
\label{sec:concrete-models-and-computations}

This section derives a formula that reconstructs a boundary torus class from relative homotopy classes associated with rectangles in a boundary annulus.
After proving the general assembly formula in Proposition~\ref{prop:assembly-rule-for-boundary-torus-class}, we apply it to the two boundary components of a braid exterior and to a knot exterior using a suitably normalized nonorientable spanning surface.

\subsection{An Assembly Formula from Local Rectangles}

The assembly formula in this subsection glues local rectangles back into one boundary torus.
Fix a filtered CW complex $E$ that models a link exterior, and use it as the target below.
Write $E^{(k)}$ for its $k$-skeleton.
In later applications, $E$ will be the chosen filtered CW model of either the braid exterior $X_P$ or the knot exterior $E_K$.
Cut the standard torus $T=S^1_\mu\times S^1_\lambda$ along the longitude $\{1\}\times S^1_\lambda$, and denote the resulting annulus by
\[
A:=S^1_\lambda\times[0,1].
\]
For a point $(z,r)$ of $A$, the coordinate $z\in S^1_\lambda$ runs in the longitude direction, and $r\in[0,1]$ runs transversely in the meridian direction.
Define the inclusions of the lower and upper boundary circles by
\[
\iota^\downarrow,\iota^\uparrow\colon S^1_\lambda\longrightarrow A,
\qquad
\iota^\downarrow(z):=(z,0),
\quad
\iota^\uparrow(z):=(z,1).
\]
Choose a CW decomposition of $A$ in which its two boundary circles are $1$-subcomplexes, and fix a filtration-preserving cellular map
\[
F_A\colon A\longrightarrow E
\]
that represents the boundary annulus in the specific model under consideration.
Then $F_A\circ\iota^\downarrow$ and $F_A\circ\iota^\uparrow$ are cellular loops in $E^{(1)}$.
We do not assume that these two loops agree pointwise.
After fixing corresponding subdivision points and basepoint paths, we assume that the loops represent the same based homotopy class.
For this reason, we do not yet take the quotient that identifies the two boundary circles of $A$.
We next decompose $A$ into rectangles arranged around the annulus and subdivide the upper and lower boundaries accordingly.
Assume that $F_A$ sends corresponding subdivision points on the upper and lower boundaries to the same $0$-cell.
We read the relative homotopy classes of each edge and each rectangle using the image of its corner as the local basepoint.
When changing to another basepoint, we use paths in a fixed maximal tree in $E^{(1)}$ and use the same path for corresponding subdivision points on the upper and lower boundaries.

We now associate local elements with the rectangles in the decomposition of $A$.
Decompose $A$ into rectangles $R_1,\ldots,R_n$ arranged in this order around the annulus.
By the assumption in the preceding paragraph, $F_A$ maps the two endpoints of the left vertical side of $R_i$ to the same vertex.
Denote this vertex by $z_i$, and denote by $m_i$ the loop obtained by applying $F_A$ to the vertical side.
Similarly, define $z_{i+1}$ and $m_{i+1}$ using the right side of $R_i$.
Denote by $u_i^\uparrow$ and $u_i^\downarrow$ the paths obtained by reading the upper and lower sides, respectively, from left to right.
The orientations of this notation are shown in the following diagram.
\begin{center}
\small
\(
\begin{array}{ccc}
z_i & \xrightarrow{\ u_i^\uparrow\ } & z_{i+1}\\[-1pt]
{\scriptstyle m_i}\Big\uparrow && \Big\uparrow{\scriptstyle m_{i+1}}\\[-1pt]
z_i & \xrightarrow{\ u_i^\downarrow\ } & z_{i+1}
\end{array}
\)
\end{center}
More precisely,
\[
m_i\in(\Pi E)_1(z_i),
\quad
u_i^\uparrow,u_i^\downarrow\in(\Pi E)_1(z_i,z_{i+1}),
\quad
\widetilde\Omega_i\in(\Pi E)_2(z_i).
\]
Here $\widetilde\Omega_i$ is the relative homotopy class in $(\Pi E)_2(z_i)$ obtained by basing the restriction $F_A|_{R_i}$ at $z_i$ along the basepoint path (whisker) chosen in the preceding paragraph.
We read the indices cyclically and set $z_{n+1}=z_1$ and $m_{n+1}=m_1$.
Orient $\widetilde\Omega_i$ so that the boundary of the rectangle is traversed along the left side, the upper side, the reverse of the right side, and the reverse of the lower side, in that order.
Then
\begin{equation}
\label{eq:local-strip-boundary}
\partial_2(\widetilde\Omega_i)
=m_i u_i^\uparrow m_{i+1}^{-1}(u_i^\downarrow)^{-1}.
\end{equation}

We next recover a torus by using a homotopy between the two boundary loops of the annulus.
Choose $1$ as the basepoint of $S^1_\lambda$, and assume that the corresponding subdivision points chosen above satisfy
\[
F_A\bigl(\iota^\uparrow(1)\bigr)
=F_A\bigl(\iota^\downarrow(1)\bigr)
=z_1.
\]
Write the upper and lower boundary restrictions as based cellular loops
\[
\gamma^\uparrow,\gamma^\downarrow
\colon(S^1_\lambda,1)\longrightarrow(E^{(1)},z_1),
\qquad
\gamma^\uparrow:=F_A\circ\iota^\uparrow,
\quad
\gamma^\downarrow:=F_A\circ\iota^\downarrow.
\]
Assume that these two loops represent the same element of $(\Pi E)_1(z_1)=\pi_1(E^{(1)},z_1)$; that is,
\begin{equation}
\label{eq:assembled-longitude}
w
:=[\gamma^\uparrow]
=u_1^\uparrow\cdots u_n^\uparrow
=u_1^\downarrow\cdots u_n^\downarrow
=[\gamma^\downarrow].
\end{equation}
This assumption requires equality of the two complete boundary words, not the equalities $u_i^\uparrow=u_i^\downarrow$ for every $i$.
Since $E^{(1)}$ is a graph, we may choose a based cellular homotopy between the two loops,
\[
H\colon
\bigl(S^1_\lambda\times[0,1],\{1\}\times[0,1]\bigr)
\longrightarrow(E^{(1)},z_1),
\qquad
H(-,0)=\gamma^\downarrow,
\quad
H(-,1)=\gamma^\uparrow.
\]
Let $C:=S^1_\lambda\times[0,1]$ be an auxiliary annulus that carries this homotopy, and denote the inclusions of its lower and upper boundaries by
\[
j_0,j_1\colon S^1_\lambda\longrightarrow C,
\qquad
j_0(z):=(z,0),
\quad
j_1(z):=(z,1).
\]
Glue the lower boundary of $C$ to the lower boundary of $A$ and set
\[
A_H:=A\cup_{\iota^\downarrow,j_0}C.
\]
Give $C$ a product CW structure for which $H$ is cellular, and regard $A$ and $C$ as subspaces of $A_H$ through their standard inclusions.
Define a cellular map $F_H\colon A_H\to E$ by
\[
F_H|_A:=F_A,
\qquad
F_H|_C:=H.
\]
The equality $H(-,0)=\gamma^\downarrow$ shows that $F_H$ is well defined on the gluing locus, and
\[
F_H\circ\iota^\uparrow
=\gamma^\uparrow
=F_H\circ j_1.
\]
Identify the upper boundary of $A$ pointwise with the upper boundary of $C$, and define
\[
T_H:=A_H/\bigl(\iota^\uparrow(z)\sim j_1(z)\bigr).
\]
Denote the quotient map by
\[
q_H\colon A_H\longrightarrow T_H.
\]
Since $F_H$ agrees pointwise on the two identified boundary components, there is a unique cellular map $f_H\colon T_H\to E$ such that
\[
F_H=f_H\circ q_H.
\]
The space $T_H$ is a torus.
We identify it with the standard torus $T=S^1_\mu\times S^1_\lambda$ by taking the orientation along the boundary circles as the $\lambda$-direction and the orientation transverse to the two annuli as the $\mu$-direction.
The $\mu$-direction first traverses a vertical side of $A$ from bottom to top in the diagram above and then returns along the basepoint line of $C$ from $j_1(1)$ to $j_0(1)$.
Finally, set $m:=m_1$, and define the partial products along the lower boundary and the changes of basepoint of the local classes by
\[
\ell_0:=1,
\qquad
\ell_i:=\ell_{i-1}u_i^\downarrow,
\qquad
\Omega_i:=\ell_{i-1}\cdot\widetilde\Omega_i,
\qquad
\Omega:=\Omega_1\cdots\Omega_n.
\]

\begin{Proposition}
\label{prop:assembly-rule-for-boundary-torus-class}
The following assembly formula computes the triple induced by the annulus.
If the resulting torus map represents the fixed boundary parametrization, this triple is the corresponding boundary torus class.
Fix, as above, an annulus map $F_A\colon A\to E$, a rectangle decomposition $A=R_1\cup\cdots\cup R_n$, local data satisfying \eqref{eq:local-strip-boundary}, and a based homotopy $H$ chosen as in the construction following \eqref{eq:assembled-longitude}.
Let $m=m_1$, let $w$ be the common boundary word, and let $\Omega$ be defined using the partial products from the initial vertex.
Let $f_H\colon T_H\to E$ be obtained by attaching to $A$ the auxiliary annulus $C$ that carries $H$.
Then
\[
(f_H)_*(\Theta)=\Omega,
\qquad
(f_H)_*(\mu)=m,
\qquad
(f_H)_*(\lambda)=w,
\]
and
\[
\partial_2(\Omega)=[m,w].
\]
If $f_H$ is homotopic to the fixed boundary parametrization, then $(\Omega,m,w)$ is the corresponding boundary torus class.

\end{Proposition}

\begin{proof}
We prove the assembly formula by multiplying the transported local classes.
The $(\Pi E)_1$-equivariance of the boundary map $\partial_2$ gives
\[
\partial_2(\Omega_i)
=\ell_{i-1}m_i u_i^\uparrow m_{i+1}^{-1}\ell_i^{-1}.
\]
Multiplying these expressions in the order corresponding to $R_1,\ldots,R_n$ cancels the intermediate terms and gives
\[
\partial_2(\Omega)
=m(u_1^\uparrow\cdots u_n^\uparrow)m^{-1}\ell_n^{-1}
=[m,w].
\]
Since the image of $H$ lies in $E^{(1)}$, the auxiliary annulus $C$ adds no relative class in degree $2$.
Applying the homotopy addition formula \cite[Section~9.9]{BHS11} to the rectangle decomposition of $A$ and the auxiliary annulus $C$ therefore gives $(f_H)_*(\Theta)=\Omega$.
Since $H$ is constant at the basepoint, crossing $C$ does not change the image of the meridian, and hence $(f_H)_*(\mu)=m$.
The two boundary circles that are glued together represent the common boundary word $w$, so $(f_H)_*(\lambda)=w$.
The final assertion follows from the definition of a boundary torus class.
\end{proof}

\subsection{Boundary Torus Classes for Braid Exteriors}
\label{subsec:boundary-torus-classes-for-braids-in-a-solid-torus}

This subsection computes the boundary torus classes of the two boundary components of a braid exterior.
We first construct the free crossed complex of the pattern exterior $X_P$ from a braid $\beta$.
We then compute the boundary torus classes of the outer boundary $\partial V$ and the inner boundary $\partial\nu P$, and correct the blackboard framing on the inner boundary to the preferred longitude.

We begin by fixing the circle direction of the mapping torus and a free basis for its fiber.
Let $V\cong S^1\times D^2$ be an oriented solid torus, and let $P=\widehat\beta\subset\operatorname{int}(V)$ be the closure of a braid $\beta\in B_n$ with $n$ strands.
Assume that $P$ is a knot, and set
\[
X_P:=V\setminus\operatorname{int}(\nu P).
\]
Regard $V$ as the solid torus obtained from $D^2\times[0,1]$ by identifying $D^2\times\{0\}$ with $D^2\times\{1\}$, and denote the interval coordinate by $s\in[0,1]$.
After this identification, call the direction in which $s$ increases the circle direction of the mapping torus.
Orient every strand of the braid, and hence its closure $P$, in this direction.
To measure the framing of $P$, embed $V$ as an unknotted solid torus in $S^3$ with its product longitude equal to the preferred longitude of its core.
Orient $D^2$ so that its orientation followed by the direction in which $s$ increases gives the orientation of $V$.
Choose pairwise disjoint closed disks $D_1,\ldots,D_n$ in $\operatorname{int}(D^2)$, each with the orientation induced from $D^2$.
Choose a basepoint $*$ on the outer boundary $\partial D^2$ of
\[
\Sigma_n:=D^2\setminus\operatorname{int}(D_1\sqcup\cdots\sqcup D_n).
\]
For each $i$, choose a point $q_i\in\partial D_i$ and fix a path
\[
a_i\colon[0,1]\longrightarrow\Sigma_n
\]
such that $a_i(0)=*$ and $a_i(1)=q_i$.
Choose these paths to be pairwise disjoint except at the basepoint and so that the holes occur in the order $D_1,\ldots,D_n$ along the orientation of $\partial D^2$.
Let $x_i$ be the loop that follows $a_i$ to $q_i$, traverses $\partial D_i$ once in the positive direction, and returns along $a_i$ in reverse.
These loops form a free basis for the fundamental group of the fiber $\Sigma_n$.
Write
\[
\pi_1(\Sigma_n,*)=\langle x_1,\ldots,x_n\rangle,
\qquad
\delta_0:=x_1\cdots x_n.
\]
Thus, each $x_i$ is the positive meridian of the $i$th hole, or equivalently, of the corresponding braid strand.

The Artin action gives a mapping-torus description of the braid exterior.
The positive Artin generator $\sigma_i$ of $B_n$ acts on the free group $\pi_1(\Sigma_n,*)=\langle x_1,\ldots,x_n\rangle$ by
\[
\sigma_i(x_i)=x_ix_{i+1}x_i^{-1},
\qquad
\sigma_i(x_{i+1})=x_i,
\qquad
\sigma_i(x_j)=x_j\quad(j\neq i,i+1).
\]
Choose a homeomorphism $h_\beta\colon\Sigma_n\to\Sigma_n$ that induces the automorphism associated with $\beta$ under this Artin action and fixes the outer boundary of $\Sigma_n$ pointwise, and write $\beta_*:=(h_\beta)_*$.
Then $\beta_*(\delta_0)=\delta_0$.
Let $\operatorname{wr}(\beta)$ be the exponent sum of any braid word representing $\beta$ in the Artin generators $\sigma_1,\ldots,\sigma_{n-1}$.
The exterior $X_P$ is the mapping torus of $h_\beta$:
\[
X_P\cong\Sigma_n\times[0,1]/(z,1)\sim(h_\beta(z),0).
\]
Denote the quotient map in this expression by
\[
\pi_\beta\colon\Sigma_n\times[0,1]\longrightarrow X_P,
\]
and take $x_0:=\pi_\beta(*,0)$ as the basepoint of $X_P$.
Define a path from $[0,1]$ to $X_P$ by
\[
\tau\colon[0,1]\longrightarrow X_P,
\qquad
\tau(s):=\pi_\beta(*,s)=[*,s].
\]
Here the brackets in $[*,s]$ denote the equivalence class of $(*,s)\in\Sigma_n\times[0,1]$ in the quotient space $X_P$.
Since $h_\beta(*)=*$,
\[
\tau(1)=\pi_\beta(*,1)
=\pi_\beta(h_\beta(*),0)
=x_0
=\tau(0).
\]
Thus, $\tau$ is a loop based at $x_0$.
Use the image of $\tau$ as the $1$-cell in the circle direction of the mapping-torus CW structure, and denote its edge-path class by $[\tau]\in(\Pi X_P)_1(x_0)$.
Define the degree-$1$ generator with the opposite orientation by
\[
t:=[\tau]^{-1}\in(\Pi X_P)_1(x_0).
\]
The brackets in $[\tau]$ denote the based edge-path class of the loop and thus have a different role from the brackets in $[*,s]$, which denote a point.
With this orientation convention, the relator associated with the $i$th fiber generator is
\[
t x_i t^{-1}\beta_*(x_i)^{-1}.
\]

\begin{Proposition}
\label{prop:braid-pattern-exterior-crossed-complex}
The crossed complex of the braid exterior has the following explicit free presentation.
The crossed complex $\Pi X_P$ is reduced and free, with degree-$1$ generators $t,x_1,\ldots,x_n$ and degree-$2$ generators $\rho_1,\ldots,\rho_n$, and
\begin{equation}
\label{eq:braid-pattern-relators}
\partial_2(\rho_i)
=t x_i t^{-1}\beta_*(x_i)^{-1}
\qquad(1\leq i\leq n).
\end{equation}

\end{Proposition}

\begin{proof}
We obtain this presentation from a cellular decomposition of the mapping torus.
Decompose $\Sigma_n$ into one $0$-handle and $n$ $1$-handles corresponding to the loops $x_1,\ldots,x_n$.
Take the product of this decomposition with $[0,1]$, and glue its $s=1$ end to its $s=0$ end by $h_\beta$.
The $1$-cell obtained by following the trajectory of the basepoint $*$ in the direction in which $s$ decreases gives the degree-$1$ generator $t$ defined above.
The trajectory in the $s$-direction of the core $1$-cell of the $i$th $1$-handle of $\Sigma_n$ gives a $2$-cell; denote the corresponding degree-$2$ generator by $\rho_i$.
Its attaching loop successively follows $t$, $x_i$, $t^{-1}$, and $\beta_*(x_i)^{-1}$ from the basepoint.
Hence
\[
\partial_2(\rho_i)=t x_i t^{-1}\beta_*(x_i)^{-1},
\]
which gives \eqref{eq:braid-pattern-relators}.
The assertion now follows from Proposition~\ref{prop:handle-classes-free-basis}.

\end{proof}

\begin{Remark}
\label{rem:braid-local-characteristic-classes}
The relative homotopy class represented by each rectangle is uniquely determined by its boundary word.
The surface $\Sigma_n$ is aspherical, so its universal cover is contractible.
The universal cover of $X_P$ is homeomorphic to $\widetilde\Sigma_n\times\mathbb R$, and hence $\pi_2(X_P)=0$.
Denote by $K_{X_P}$ the $2$-dimensional CW model obtained in the proof of Proposition~\ref{prop:braid-pattern-exterior-crossed-complex}.
This model satisfies $K_{X_P}=K_{X_P}^{(2)}$ and $K_{X_P}\simeq X_P$.
Thus, the exact sequence of the pair $(K_{X_P},K_{X_P}^{(1)})$, together with $\pi_2(X_P)=0$, shows that
\[
\partial_2\colon(\Pi X_P)_2\longrightarrow(\Pi X_P)_1
\]
is injective.
This injectivity implies that the word obtained by reading the oriented boundary of a rectangle uniquely determines the relative homotopy class represented by its restriction map.

\end{Remark}

We now specify the local elements for the rectangles in the outer boundary.
For each rectangle, set
\[
m_i:=t,
\qquad
u_i^\uparrow:=x_i,
\qquad
u_i^\downarrow:=\beta_*(x_i),
\qquad
\widetilde\Omega_i^A:=\rho_i.
\]
The identity $\partial_2(\rho_i)=t x_i t^{-1}\beta_*(x_i)^{-1}$ agrees with \eqref{eq:local-strip-boundary} after substituting $m_i=m_{i+1}=t$, $u_i^\uparrow=x_i$, and $u_i^\downarrow=\beta_*(x_i)$.
Choose a $1$-dimensional spine of $\Sigma_n$ so that the word representing its outer boundary is $x_1\cdots x_n$.
For the $i$th rectangle in the outer boundary, the word read along its four sides also agrees with the right-hand side of \eqref{eq:braid-pattern-relators}.
The uniqueness in Remark~\ref{rem:braid-local-characteristic-classes} therefore shows that the relative homotopy class represented by this rectangle is $\rho_i=\widetilde\Omega_i^A$.
The upper boundary word is $\delta_0=x_1\cdots x_n$, and the lower boundary word is $\beta_*(\delta_0)=\delta_0$.
Thus, their common boundary word is $\delta_0$.
Although the individual edge words $x_i$ and $\beta_*(x_i)$ need not agree, the assembly formula requires only equality of the complete boundary words.
These rectangles form an annulus obtained by mapping the outer boundary cellularly into the mapping-torus CW model.
Attaching the auxiliary annulus that represents the based homotopy in \eqref{eq:assembled-longitude} gives a torus map $f_H$.
Since $h_\beta$ fixes the outer boundary of $\Sigma_n$ pointwise, this map is homotopic to the outer boundary parametrization before exchanging the coordinates.
Let $\ell_i^A$ be the partial product obtained by following the lower boundary from the first edge through the $i$th edge.
Define
\[
\ell_0^A:=1,
\qquad
\ell_i^A:=\ell_{i-1}^A\beta_*(x_i),
\qquad
\Omega_A
:=\prod_{i=1}^n
\bigl(\ell_{i-1}^A\cdot\widetilde\Omega_i^A\bigr).
\]

To obtain the coordinates used in the satellite construction, we must exchange the two coordinates produced by the outer boundary annulus.
Applying the assembly formula to the annulus constructed above gives a map $f_H\colon T\to X_P$ that sends the two generators of the standard torus to
\[
(f_H)_*(\mu)=t,
\qquad
(f_H)_*(\lambda)=\delta_0.
\]
The satellite construction instead uses the outer boundary coordinates in the order
\[
(\mu_V,\lambda_V):=(\delta_0,t).
\]
We therefore exchange the two coordinates defined by $f_H$.
For this purpose, consider the map on the standard torus $T=S^1_\mu\times S^1_\lambda$ given by
\[
J\colon T\longrightarrow T,
\qquad
J\bigl(e^{2\pi \mathrm{i}u},e^{2\pi \mathrm{i}v}\bigr)
:=\bigl(e^{2\pi \mathrm{i}v},e^{2\pi \mathrm{i}u}\bigr).
\]
Since $J$ exchanges the two circle factors,
\[
J_*(\mu)=\lambda,
\qquad
J_*(\lambda)=\mu,
\qquad
\deg J=-1.
\]
Moreover,
\[
\partial_2\bigl(J_*(\Theta)\bigr)
=[\lambda,\mu]
=[\mu,\lambda]^{-1}
=\partial_2(\Theta^{-1}).
\]
Since $\pi_2(T)=0$, the map $\partial_2\colon(\Pi T)_2\to(\Pi T)_1$ is injective, and hence $J_*(\Theta)=\Theta^{-1}$.
It follows that the map $f_{\mathrm{out}}:=f_H\circ J$ obtained after exchanging the coordinates defines the triple
\[
\Omega_{\mathrm{out}}:=\Omega_A^{-1},
\qquad
\meri_{\mathrm{out}}:=\delta_0,
\qquad
\longi_{\mathrm{out}}:=t.
\]
If $f_{\mathrm{out}}$ is taken as the outer boundary parametrization, its degree with respect to the orientation of the outer boundary is $-1$.
On the other hand, the usual basis compatible with the orientation of the solid torus $V$ is $(\delta_0,t^{-1})$.
This difference occurs because we defined $t=[\tau]^{-1}$ when fixing the circle direction of the mapping torus, so $\lambda_V=t$ is oriented opposite to the direction in which the interval coordinate $s$ increases.

\begin{Proposition}
\label{prop:outer-braid-datum}
With respect to the chosen basis $(\mu_V,\lambda_V)$, the boundary torus class of the outer boundary is $(\Omega_{\mathrm{out}},\meri_{\mathrm{out}},\longi_{\mathrm{out}})$.

\end{Proposition}

\begin{proof}
The result follows by exchanging the two coordinates in the triple obtained before the coordinate change.
Applying Proposition~\ref{prop:assembly-rule-for-boundary-torus-class} to the annulus data before exchanging coordinates gives the triple $(\Omega_A,t,\delta_0)$.
Exchanging the two circle coordinates transforms this triple as follows:
\[
(\Omega_A,t,\delta_0)
\longmapsto
(\Omega_A^{-1},\delta_0,t).
\]
This is the asserted triple.

\end{proof}

For the inner boundary, we first normalize the conjugating words of the Artin action.
Write $\bar\beta$ for the permutation of the holes and $\varepsilon_j$ for the exponent sum of $x_j$ in the free group $\pi_1(\Sigma_n,*)$.
Choose $u_i(\beta)$ by
\[
\beta_*(x_i)=u_i(\beta)x_{\bar\beta(i)}u_i(\beta)^{-1},
\qquad
\varepsilon_{\bar\beta(i)}(u_i(\beta))=0.
\]
The centralizer of $x_{\bar\beta(i)}$ in the fiber group is $\langle x_{\bar\beta(i)}\rangle$.
Thus all conjugators differ by right multiplication by a power of $x_{\bar\beta(i)}$, and exactly one satisfies the exponent condition.

Set $\varepsilon=\sum_{j=1}^n\varepsilon_j$.
For braids $\phi,\psi$, the normalized words satisfy
\[
u_i(\phi\psi)=\phi_*(u_i(\psi))\,u_{\bar\psi(i)}(\phi)
\qquad\bigl((\phi\psi)_*=\phi_*\circ\psi_*\bigr).
\]
Indeed, the right-hand side conjugates $x_{\bar\phi(\bar\psi(i))}$ to $(\phi\psi)_*(x_i)$ and has zero exponent sum in that letter, since $\phi_*$ permutes the meridians in the abelianization.
Consequently $\beta\mapsto\sum_i\varepsilon(u_i(\beta))$ is a homomorphism $B_n\to\Z$.
For the positive generator $\sigma_j$, its only nontrivial normalized conjugator is $u_j(\sigma_j)=x_j$.
It follows that
\[
\sum_{i=1}^n\varepsilon(u_i(\beta))=\operatorname{wr}(\beta).
\]

Since $P$ is a knot, $\bar\beta$ is an $n$-cycle.
Choose $i_1$, put $i_{k+1}=\bar\beta(i_k)$, and set $v_k=t^{-1}u_{i_k}(\beta)$ for $1\le k\le n$.
In particular $i_{n+1}=i_1$.
Equation~\eqref{eq:braid-pattern-relators} gives
\[
\partial_2(t^{-1}\cdot\rho_{i_k})=x_{i_k}v_kx_{i_{k+1}}^{-1}v_k^{-1}.
\]
For the $k$th rectangle choose a representative of $t^{-1}\cdot\rho_{i_k}$ with vertical edges $x_{i_k},x_{i_{k+1}}$ and both horizontal edges $v_k$.
The boundary formula permits these edge paths, and Remark~\ref{rem:braid-local-characteristic-classes} gives uniqueness of the relative class.
Glue the rectangles cyclically and apply Proposition~\ref{prop:assembly-rule-for-boundary-torus-class}, with
\[
m_k=x_{i_k},\qquad
u_k^\uparrow=u_k^\downarrow=v_k,\qquad
\widetilde\Omega_k=t^{-1}\cdot\rho_{i_k}.
\]
Denote the resulting degree-$2$ element by $\Omega_{\mathrm{bb}}$ and set
\[
\meri_{\mathrm{in}}:=x_{i_1},\qquad
w_{\mathrm{bb}}:=v_1\cdots v_n .
\]

We identify the framing of this torus map in $\pi_1(X_P,x_0)$.
Use the path $a_{i_1}$ to base the preferred longitude $\lambda_{\mathrm{pref}}$ of $P$, oriented as $P$, at $x_0$.
The relations $x_{i_k}v_k=v_kx_{i_{k+1}}$ show that $w_{\mathrm{bb}}$ commutes with $\meri_{\mathrm{in}}$.
Both $w_{\mathrm{bb}}$ and $\lambda_{\mathrm{pref}}$ project to $n$ times the positive generator of the mapping-torus circle.
Hence $w_{\mathrm{bb}}\lambda_{\mathrm{pref}}^{-1}$ lies in the fiber group and centralizes $x_{i_1}$, so it is a power of $\meri_{\mathrm{in}}$.
To compute the exponent, take linking numbers with $P$ in the chosen embedding of $V$.
Each $x_i$ has linking number $1$, while $t$ has linking number $0$ because the product longitude of $\partial V$ bounds a disk outside $V$.
The exponent is therefore $\sum_i\varepsilon(u_i(\beta))=\operatorname{wr}(\beta)$.
Thus $w_{\mathrm{bb}}$ represents the blackboard longitude of the standard closed braid diagram \cite{Lic97}.
The assembled torus map and the blackboard boundary parametrization now induce the same based map on $\pi_1$.
Since $\pi_2(X_P)=0$, these maps are homotopic with the chosen basepoint path, so $(\Omega_{\mathrm{bb}},\meri_{\mathrm{in}},w_{\mathrm{bb}})$ is the boundary torus class for this framing.

Put $e:=-\operatorname{wr}(\beta)$ and define
\[
\longi_{\mathrm{in}}:=\meri_{\mathrm{in}}^{\,e}w_{\mathrm{bb}},
\qquad
\Omega_{\mathrm{in}}:=\meri_{\mathrm{in}}^{\,e}\cdot\Omega_{\mathrm{bb}}.
\]

\begin{Proposition}
\label{prop:inner-braid-datum}
With respect to the preferred parametrization, the boundary torus class of the inner boundary $\partial\nu P$ is $(\Omega_{\mathrm{in}},\meri_{\mathrm{in}},\longi_{\mathrm{in}})$.
Here $\meri_{\mathrm{in}}$ is the positive meridian, and $\longi_{\mathrm{in}}$ is the preferred longitude oriented as $P$.
The degree of this parametrization is $-1$ with respect to the orientation of the inner boundary induced from $X_P$.

\end{Proposition}

\begin{proof}
The meridional twist with exponent $e$ changes the blackboard longitude to the preferred longitude. Applying Remark~\ref{rem:meridional-dehn-twist} to $(\Omega_{\mathrm{bb}},\meri_{\mathrm{in}},w_{\mathrm{bb}})$ gives the asserted triple.
The ordered positive meridian and longitude oriented as $P$ give the boundary orientation of $\nu P$, opposite to that induced from $X_P$.

\end{proof}

\subsection{Boundary Torus Classes from Disk-and-Band Decompositions in Crosscap Normal Form}
\label{subsec:the-non-orientable-seifert-surface-model}

This subsection computes a boundary torus class from a spanning surface equipped with a disk-and-band decomposition in crosscap normal form.
We first construct the free crossed complex in degrees $1$ and $2$ of the knot exterior from this decomposition, and then determine the local class represented by each rectangle in the boundary annulus.
We combine these classes by the assembly formula and finally correct the surface framing to the preferred longitude.

We use the following disk-and-band decomposition in crosscap normal form.
Let $F\subset S^3$ be a nonorientable surface that spans a knot $K$, and fix an orientation of $K=\partial F$.
We call $F$ admissible if $\overline{S^3\setminus\nu F}$ is a handlebody, and assume throughout this subsection that $F$ is admissible.
Assume that $F$ has a disk-and-band decomposition
\[
F=D^2\cup B_1\cup\cdots\cup B_b,
\qquad
b=\beta_1(F),
\]
that satisfies the following two conditions.
\begin{enumerate}[(1)]
\item For each $B_i$, close its core arc by an arc in $D^2$ that joins the two endpoints.
The resulting simple closed curve is one-sided in $F$.
We call $B_i$ a one-sided band in this sense.
\item In the cyclic order on $\partial D^2$, the attaching intervals occur in the order
$B_1,B_1,B_2,B_2,\ldots,B_b,B_b$.
Thus, the attaching intervals of distinct bands do not interlace, and the two free boundary arcs of each band form one consecutive block along $\partial F$.
\end{enumerate}
We call a decomposition satisfying these two conditions a disk-and-band decomposition in crosscap normal form.

The regular neighborhood of $F$ determines a handle decomposition of the knot exterior.
Let $\nu F$ be a closed regular neighborhood of $F$, and set
\[
E_K:=S^3\setminus\operatorname{int}(\nu K),
\qquad
\mathcal H_2:=\nu F,
\qquad
\mathcal H_1:=\overline{S^3\setminus\nu F}.
\]
Let $I=[-1,1]$.
Thickening the disk-and-band decomposition of $F$ to a regular neighborhood makes $D^2\times I$ a $0$-handle of $\mathcal H_2$ and makes each $B_i\times I$ a $1$-handle of $\mathcal H_2$.
Since $F$ is nonorientable, $\mathcal H_2$ is globally the nontrivial $I$-bundle over $F$, not the product $F\times I$.
The product descriptions above remain valid on the individual handles.
Thus, $\mathcal H_2$ is a handlebody of genus $b$, and the admissibility assumption implies that $\mathcal H_1$ is also a handlebody of genus $b$.
Together they give a Heegaard splitting $S^3=\mathcal H_1\cup\mathcal H_2$ adapted to $F$.
Choose a handle decomposition of $\mathcal H_1$ consisting of one $0$-handle and $b$ $1$-handles.

The handle decomposition of $E_K$ is obtained by reading this decomposition dually from the $\mathcal H_1$ side.
The subspace $(D^2\times I)\cap E_K$ is attached to $\mathcal H_1$ along the two disk faces $D^2\times\{-1\}$ and $D^2\times\{1\}$ and becomes a $1$-handle of $E_K$.
Denote its degree-$1$ class by $\meri_F$, with the orientation chosen to agree with the positive meridian determined by the orientation of $K$.
Denote the degree-$1$ classes determined by the chosen $b$ $1$-handles of $\mathcal H_1$ by $g_1,\ldots,g_b$.
The part arising from each $B_i\times I$ becomes a $2$-handle of $E_K$; denote its degree-$2$ class by $r_i$.
These handles give all of $E_K$.
The $0$-handle part of $\mathcal H_2$ already appears as the $1$-handle corresponding to $\meri_F$, so no $3$-handle remains.
Let $K_F$ be the CW model associated with this handle decomposition, and let
\[
h_F\colon K_F\longrightarrow E_K
\]
be the specified homotopy equivalence.
For the remainder of this subsection, write $\Pi E_K:=\Pi K_F$.

We cut the boundary torus along the surface push-off and read the two sides of each band in the fixed CW model.
Put $F'=F\cap E_K$ and orient its boundary as $K$. This push-off determines the surface framing.
The horizontal boundary of a regular $I$-bundle neighborhood $N$ of $F'$ is the orientation double cover of $F'$.
The two lifts of $\partial F'$, oriented as $K$, bound the complementary annulus $\overline{\partial E_K\setminus N}$ on the $\mathcal H_1$ side.
Here we identify $\overline{E_K\setminus N}$ with $\mathcal H_1$ by adjusting a collar near $K$.
Use the subarc of a meridian in this annulus to identify the two lifts of the basepoint, and choose the whiskers in the handle model compatibly with this arc.
Let $A_F$ be the annulus obtained by cutting $\partial\nu K$ along the surface push-off $\partial F'$.
After assigning the parts arising from $D^2$ to adjacent rectangles, the crosscap normal form gives a cyclic decomposition
\[
A_F=Q_1\cup\cdots\cup Q_b.
\]
Here $Q_i$ follows the two free boundary arcs of $B_i$ and passes once through each face of $B_i\times I$.
Through the fixed homotopy equivalence $h_F\colon K_F\to E_K$, choose a cellular map $A_F\to K_F$ that represents the regluing map and sends the parts arising from $D^2$ into $K_F^{(1)}$.
Write $x_F$ for the unique vertex of $K_F$.
Map the two sides of $Q_i$ that run along $B_i$ into $K_F^{(1)}$, base them at $x_F$ using the fixed whiskers, and denote the resulting words by
\[
u_i^\sharp,u_i^\flat\in\langle g_1,\ldots,g_b\rangle.
\]
Choose their orientations so that the attaching loop of $r_i$ reads $\meri_Fu_i^\sharp\meri_F(u_i^\flat)^{-1}$.
Since $B_i$ is one-sided, both meridional segments occur with positive orientation; hence
\begin{equation}
\label{eq:crosscap-relators}
\partial_2(r_i)
=\meri_Fu_i^\sharp\meri_F(u_i^\flat)^{-1}.
\end{equation}
With these basepoints, the complementary annulus gives a based homotopy between the two boundary lifts.
Their complete words therefore agree in the free group $\pi_1(\mathcal H_1)=\langle g_1,\ldots,g_b\rangle$, as required by the assembly formula.

\begin{Proposition}
\label{prop:nonorientable-surface-degree12-model}
The crossed complex $\Pi E_K$ is reduced and free.

Its degree-$1$ generators are $\meri_F,g_1,\ldots,g_b$.
Its degree-$2$ generators are $r_1,\ldots,r_b$, with boundaries given in \eqref{eq:crosscap-relators}.

\end{Proposition}

\begin{proof}
This presentation follows from the handle decomposition in crosscap normal form.
The dual decomposition of the Heegaard splitting adapted to $F$ has one $0$-handle, the $1$-handles corresponding to $\meri_F,g_1,\ldots,g_b$, the $2$-handles corresponding to the bands, and no $3$-handle.
The attaching loop of each band successively follows $\meri_F$, $u_i^\sharp$, $\meri_F$, and $(u_i^\flat)^{-1}$.
Since the band is one-sided, both meridians occur with positive orientation, which gives \eqref{eq:crosscap-relators}.
The assertion follows from Proposition~\ref{prop:handle-classes-free-basis}.

\end{proof}

Each local rectangle determines an explicit element whose boundary has the form required by the assembly formula.
The rectangle $Q_i$ passes through one face of the band as $r_i$ and through the opposite face as $r_i^{-1}$ transported along $u_i^\flat\meri_F^{-1}$.
Following this geometric order, define
\[
\widetilde\Omega_{F,i}
:=
r_i\bigl(u_i^\flat\meri_F^{-1}\cdot r_i^{-1}\bigr),
\qquad
u_{F,i}^\uparrow:=u_i^\sharp u_i^\flat,
\qquad
u_{F,i}^\downarrow:=u_i^\flat u_i^\sharp,
\qquad
m_{F,i}:=\meri_F.
\]
Read the indices cyclically, and set $m_{F,b+1}:=m_{F,1}=\meri_F$.
Using \eqref{eq:crosscap-relators}, we compute
\[
\begin{aligned}
\partial_2(\widetilde\Omega_{F,i})
&=
\partial_2(r_i)
(u_i^\flat\meri_F^{-1})
\partial_2(r_i)^{-1}
(\meri_F(u_i^\flat)^{-1})\\
&=
\meri_Fu_{F,i}^\uparrow
\meri_F^{-1}
(u_{F,i}^\downarrow)^{-1}\\
&=
m_{F,i}u_{F,i}^\uparrow
m_{F,i+1}^{-1}
(u_{F,i}^\downarrow)^{-1}.
\end{aligned}
\]

\begin{Remark}
\label{rem:crosscap-local-class-uniqueness}
The geometric description of the local class also proves its uniqueness.
Consider the rectangle $Q_i\subset A_F$ fixed above.
Using the whisker from $x_F$ to the lower-left corner of $Q_i$, regard the restriction to $Q_i$ of the cellular map chosen above as a relative map
\[
q_{F,i}\colon(I^2,\partial I^2)\longrightarrow(K_F,K_F^{(1)}).
\]
The map $q_{F,i}$ passes through one face of the band as $r_i$ and through the opposite face as $r_i^{-1}$ transported along $u_i^\flat\meri_F^{-1}$.
Hence
\[
[q_{F,i}]
=r_i\bigl(u_i^\flat\meri_F^{-1}\cdot r_i^{-1}\bigr)
=\widetilde\Omega_{F,i}
\quad\text{in }\pi_2(K_F,K_F^{(1)},x_F).
\]
Reading the oriented boundary of $Q_i$ gives, in order, $\meri_F$, $u_{F,i}^\uparrow$, $\meri_F^{-1}$, and $(u_{F,i}^\downarrow)^{-1}$.
This word agrees with the boundary computed above.
We have $K_F=K_F^{(2)}\simeq E_K$, and the knot exterior is irreducible.
The sphere theorem therefore gives $\pi_2(K_F)=0$ \cite{Lic97}.
The exact sequence of the pair $(K_F,K_F^{(1)})$ shows that $\partial_2\colon(\Pi E_K)_2\to(\Pi E_K)_1$ is injective.
Thus, $\widetilde\Omega_{F,i}$ is the unique element in degree $2$ with this boundary word.
The rectangles $Q_1,\ldots,Q_b$ form a rectangle decomposition of the actual boundary annulus $A_F$.

\end{Remark}

The assembly formula now combines the local classes for the surface framing.
Let
\[
w_F
:=
u_{F,1}^\uparrow\cdots u_{F,b}^\uparrow
=
u_{F,1}^\downarrow\cdots u_{F,b}^\downarrow.
\]
This common word represents the push-off determined by the surface framing.
Apply the product formula of Proposition~\ref{prop:assembly-rule-for-boundary-torus-class} to the local data above, and denote its degree-$2$ output by $\Omega_F^{\mathrm{sf}}$.
The annulus $A_F$ is obtained by cutting the actual boundary torus along this push-off, and $\pi_2(K_F)=0$ by Remark~\ref{rem:crosscap-local-class-uniqueness}.
Hence the torus map reconstructed by the assembly formula is homotopic relative to the basepoint to the boundary parametrization determined by the surface framing, and
\[
(\Omega_F^{\mathrm{sf}},\meri_F,w_F)
\]
is its boundary torus class.

We correct the surface framing to the preferred longitude.
Let $\lambda_{\mathrm{pref}}$ be the preferred longitude, oriented as $K$.
Since $w_F$ also represents a push-off oriented as $K$, there is a unique $e_F\in\Z$ such that
\[
[w_F]
=e_F[\meri_F]+[\lambda_{\mathrm{pref}}]
\quad\text{in }H_1(\partial\nu K;\Z).
\]
Therefore,
\[
\longi_F:=\meri_F^{-e_F}w_F
\]
represents the preferred longitude.
Following Remark~\ref{rem:meridional-dehn-twist}, define
\[
\Omega_F:=\meri_F^{-e_F}\cdot\Omega_F^{\mathrm{sf}}.
\]

\begin{Proposition}
\label{prop:nonorientable-surface-boundary-class}
The boundary torus class for the preferred parametrization of $\partial\nu K$ is $(\Omega_F,\meri_F,\longi_F)$.

\end{Proposition}

\begin{proof}
The preferred boundary torus class follows from the class for the surface framing by a meridional Dehn twist.
Apply the Dehn twist about the meridian with exponent $-e_F$ to the triple $(\Omega_F^{\mathrm{sf}},\meri_F,w_F)$ associated with the boundary parametrization determined by the surface framing.
The transformation formula in Remark~\ref{rem:meridional-dehn-twist} gives $(\Omega_F,\meri_F,\longi_F)$.

\end{proof}

\section*{Acknowledgments}

The author thanks Takefumi Nosaka for suggesting the idea of constructing signed Christoffel words from continued fraction expansions.

\section*{Use of generative AI}

The author used OpenAI's ChatGPT in preparing this manuscript.
The author takes full responsibility for the mathematical content of the manuscript.

\end{document}